\documentclass[11pt]{amsart}
\usepackage[T1]{fontenc}
\usepackage{lmodern,amsmath,amssymb,mathtools,array,booktabs,longtable}
\usepackage{geometry}
\usepackage[expansion=false]{microtype}
\usepackage[colorlinks=true,linkcolor=blue,citecolor=blue,urlcolor=blue]{hyperref}
\hypersetup{bookmarksdepth=3}

\usepackage{mathrsfs}

\newtheorem{theorem}{Theorem}[section]
\newtheorem{proposition}[theorem]{Proposition}
\newtheorem{lemma}[theorem]{Lemma}
\newtheorem{corollary}[theorem]{Corollary}
\theoremstyle{definition}
\newtheorem{definition}[theorem]{Definition}
\newtheorem{example}[theorem]{Example}
\theoremstyle{remark}
\newtheorem{remark}[theorem]{Remark}

\newcommand{\Pp}{\mathcal{P}}
\newcommand{\eps}{\varepsilon}
\newcommand{\Id}{\operatorname{id}}
\newcommand{\cf}{\operatorname{cf}}

\newcommand{\End}{\operatorname{End}}
\newcommand{\Hom}{\operatorname{Hom}}
\newcommand{\ad}{\operatorname{ad}}

\newcommand{\Tr}{\operatorname{Tr}}
\newcommand{\rank}{\operatorname{rank}}
\newcommand{\Irr}{\operatorname{Irr}}
\newcommand{\corad}{\operatorname{corad}}
\newcommand{\supp}{\operatorname{supp}}
\newcommand{\Mat}{\operatorname{M}}
\newcommand{\Mco}{\mathcal M^c}
\newcommand{\Cl}{\operatorname{Cl}}
\newcommand{\Aut}{\operatorname{Aut}}
\newcommand{\GL}{\operatorname{GL}}
\newcommand{\gr}{\operatorname{gr}}
\newcommand{\id}{\operatorname{id}}
\newcommand{\C}{\mathcal{C}}
\newcommand{\Lie}{\operatorname{Lie}}
\newcommand{\Spec}{\operatorname{Spec}}
\newcommand{\Imax}{I_{\max}}
\newcommand{\B}{\mathcal B}
\newcommand{\YD}{{}^{k\Gamma}_{k\Gamma}\mathcal{YD}}

\title[Coradically generated Hopf algebras]{Finite-dimensional Hopf algebras generated by the coradical}
\author{Rongchuan Xiong}
\address{Department of Mathematics, Changzhou University,
	Changzhou 213164, China}
\email{rcxiong@foxmail.com}

\author{Naihong Hu}
\address{School of Mathematical Sciences, East China Normal
	University, Shanghai 200241, China}
\email{nhhu@math.ecnu.edu.cn}

\begin{document}
	\begin{abstract}
		Let $H$ be a finite-dimensional Hopf algebra over a field $k$. The paper is concerned with the question of when $H$ is generated by the coradical $H_0$. We give several characterizations of this property, involving the largest Hopf ideal contained in the Jacobson radical $J(H^*)$ and the coefficient coalgebra of the injective hull of the trivial comodule.  If the canonical map $H_0\to H/J(H)$ is surjective, then $H$ is generated by $H_0$ if and only if $H=H_0^3+J(H)^2$, where $H_0^3$ is the span of products of three elements of $H_0$.
		
		For finite locally free Hopf algebras of constant rank over a $k$-scheme in characteristic zero, including duals of liftings of a given Nichols algebra,  the locus where this property holds is Zariski open and compatible with arbitrary base change. 
		
		For a finite-dimensional pointed Hopf algebra $K$ with abelian coradical over an algebraically closed field of characteristic zero, the largest Hopf ideal contained in $J(K)$ is determined. $K^*$ is generated by its coradical if and only if  $J(K)$ contains no nonzero skew-primitive elements.  For liftings of quantum linear spaces and standard liftings of finite Cartan type,  this condition is expressed in terms of the deformation parameters of $K$.
	\end{abstract}
	
	\maketitle

	\section{Introduction}
	
	Let $H$ be a finite-dimensional Hopf algebra over a field $k$. The subalgebra $B(H)$ generated by the coradical $H_0$ is a Hopf subalgebra, called the Hopf coradical of $H$. Andruskiewitsch and Cuadra introduced the standard filtration associated with $B(H)$ and generalized the lifting method to study the structure of $H$. The fundamental problem in the generalized lifting method is to classify those generated by a fixed cosemisimple coalgebra \cite[Question~I]{AC}. Over an algebraically closed field of characteristic zero, this problem also includes the classification of finite-dimensional semisimple Hopf algebras, which remains largely open; see also \cite{AngionoSurvey}.
	
	The structure of finite-dimensional non-cosemisimple Hopf algebras generated by the coradical is still poorly understood. Over an algebraically closed field of characteristic zero, Beattie showed that a finite-dimensional bosonization $\mathcal B(V)\#k\Gamma$ with $\Gamma$ finite abelian has pointed dual \cite[Corollary 2.3]{Beattie}. For a nontrivial lifting $A$ whose coalgebra projection $A\to A_0$ is not an algebra map, $A^*$ is not pointed \cite[Proposition 2.4]{Beattie}. Radford algebras provide the most fundamental examples whose dual is generated by the coradical. More generally, explicit examples of such Hopf algebras arise mostly from basic Hopf algebras.  Andruskiewitsch and Angiono raised the question of determining when a finite-dimensional basic Hopf algebra is generated by the coradical, in terms of the deformation parameters of the dual pointed Hopf algebra \cite[Question~1]{AA}. Angiono subsequently formulated the corresponding question for families of pointed Hopf algebras \cite[Question~3.3]{AngionoSurvey}.
	
	In this paper, we give several characterizations of finite-dimensional Hopf algebras that are generated by the coradical (we call this property (CG)). We also study topological properties of the locus where (CG) holds for a family of finite-dimensional Hopf algebras parametrized by a $k$-scheme, and address \cite[Question~1]{AA} for pointed Hopf algebras with abelian coradical.

	A first characterization of (CG) is available when the canonical map $H_0\to H/J(H)$ is surjective over a field. In that case $H$ is generated by the coradical precisely when $H=H_0^3+J(H)^2$, where $H_0^3$ is the linear span of products of three elements of $H_0$, see Theorem~\ref{thm:four}. Equivalently, the multiplication map $m_3:H_0^{\otimes3}\to H/J(H)^2$ is surjective. Thus (CG) can be checked by computing the rank of $m_3$. The hypothesis on the map $H_0\to H/J(H)$ is satisfied whenever $(H^*)_0$ is a semisimple subalgebra of $H^*$; basic Hopf algebras with $\operatorname{char} k\nmid\dim k G(H^*)$ are an instance.
	
	Another characterization of (CG) uses the injective hull $E_H(k)$ of the trivial left $H$-comodule. Let $\C(H)$ be the Hopf subalgebra generated by the coefficient coalgebra $\operatorname{cf}(E_H(k))$ over any field. We show that  $H=H_0\C(H)=\C(H)H_0$ and that 
	$H$ is generated by the coradical if and only if $\C(H)\subseteq B(H)$, see Theorem~\ref{thm:injective-product}.
	In characteristic zero, we describe the largest semisimple Hopf subalgebra in terms of the square of the antipode, and analyze products of a semisimple normal Hopf subalgebra and a Hopf subalgebra with the Chevalley property.
	
	Topological properties of (CG) are studied for families of finite-dimensional Hopf algebras parametrized by schemes, including duals of liftings of a given finite-dimensional Nichols algebra. Let $k$ be a field of characteristic zero and $\mathcal H$ a finite locally free Hopf algebra of constant rank over a $k$-scheme $X$. For a point $x\in X$, write $\mathcal H_{\bar x}$ for the corresponding geometric fiber. The locus 
	$X^\circ=\{x\in X\mid B(\mathcal H_{\bar x})=\mathcal H_{\bar x}\}$
	is Zariski open, and its formation is compatible with arbitrary base change; see Theorem~\ref{thm:family}.
	
	The proof rests on the regular trace form of the dual algebra $\mathcal H^\vee$. On each fiber $(\mathcal H_{\bar x})^*$, the radical
	of this form is exactly the Jacobson radical
	$J((\mathcal H_{\bar x})^*)$. The annihilator of this
	Jacobson radical in $\mathcal H_{\bar x}$ is
	$(\mathcal H_{\bar x})_0$. Composing the resulting map with multiplication turns (CG) into a rank condition.
	
	Let $k$ be algebraically closed of characteristic zero, and $A$ a finite-dimensional pointed Hopf algebra with abelian coradical $k\Gamma$. By \cite[Theorem 4.15]{AngionoRelations}, $\gr A\cong\B(V)\#k\Gamma$, where $\B(V)$ is a Nichols algebra of diagonal type. We study (CG) for the basic Hopf algebra $A^*$. 
	
	Write $\Imax(A)$ for the largest Hopf ideal contained in $J(A)$,  studied by Chen and Hiss \cite{CH}. Under the annihilator correspondence between Hopf subalgebras and Hopf ideals, $\Imax(A)=B(A^*)^\perp$; so $A^*$ is generated by the coradical exactly when $\Imax(A)=0$.  Theorem~\ref{thm:first-layer} shows that every Hopf ideal meeting $k\Gamma$ trivially is generated by the skew-primitive elements it contains. Since $J(A)\cap k\Gamma=0$, it follows that
	$\Imax(A)$ is generated by the subspaces
	$J(A)\cap E_{g,\chi}$, where $E_{g,\chi}$ is the
	space of skew-primitive generators of type $(g,\chi)$.  $A^*$ is generated by the coradical if and only if $J(A)\cap E_{g,\chi}=0$ for all Yetter--Drinfeld pairs $(g,\chi)$.

	We next express this condition in terms of the deformation parameters. For liftings whose defining relations have deformation terms in $k\Gamma$, we construct linear contraction maps $\Psi_E$ to determine the intersections $J(A)\cap E_{g,\chi}$. When the chosen space of relations $W$ is preserved by Yetter--Drinfeld automorphisms, the maps $\Psi_E$ detect whether $J(A)\cap E_{g,\chi}$ is zero, and the criterion becomes that each $\Psi_E$ is injective, for every isotypic component $E=V_{g}^{\chi}$, see Theorem~\ref{thm:contraction}. 
	
	For liftings of quantum linear spaces and standard
	liftings of finite Cartan type, the criterion is made
	explicit in the deformation parameters: it gives
	nonvanishing and rank conditions in the first case,
	and a support condition in the second. For a fixed Cartan datum, admissible parameters satisfying the criterion are dense open when they exist. See Subsections \ref{sec:qls} and \ref{sec:cartan}.
	
	Finally, we study the coradical of the duals of liftings
	$A$ of quantum linear spaces. The coalgebra decomposition of
	such duals was studied by Beattie under additional assumptions
	\cite{Beattie}. We decompose $A$ by central idempotents.
	The coradical $(A^*)_0$ is the direct sum of the duals
	of the semisimple quotients of the resulting algebra blocks. 
	When all generators have the same character $\chi$ of order two,
	each algebra block is a Clifford algebra and its Jacobson radical and semisimple
	quotient are described explicitly. We also consider condition~(PC), which requires that for
	every generator $x_i$, there exists
	$\eta_i\in\widehat\Gamma\setminus\Omega$ such that
	$\omega(\eta_i,x_i)\in C(\bar\eta_i)_0$, where
	$\Omega\subseteq\widehat\Gamma$ parametrizes $G(A^*)$.
	Condition~(PC) implies that $A^*$ is generated by the
	coradical, but the converse does not hold, as a $48$-dimensional example shows.
	
	The paper is organized as follows: Section \ref{sec:prelim} recalls coefficient coalgebras, Yetter--Drinfeld modules and Nichols algebras, and bosonization and liftings. Section \ref{sec:structure} gives several characterizations of (CG). Section \ref{sec:families} studies the topological properties of the locus where (CG) holds for families over a scheme. Section \ref{sec:nichols} reduces the criterion to the skew-primitive elements of liftings of Nichols algebras via the contraction maps. Section \ref{sec:parameters} makes the criterion explicit in the deformation parameters for quantum linear spaces and finite Cartan type, and Section \ref{sec:blocks} studies central decompositions of the duals of liftings of quantum linear spaces and determines their coradicals in the Clifford case.
	
	\section{Preliminaries}
	\label{sec:prelim}

	\subsection*{Notation}
	
	Unless otherwise stated, the base field $k$ is algebraically closed 
	of characteristic zero. Write $C_m$ for a cyclic group of order $m$, and
	$\Mat_d(k)$ for the algebra of $d\times d$ matrices over
	$k$. The standard matrix units and the identity matrix
	are denoted by $e_{ij}$ and $I_d$, respectively.

	Let $A$ be an algebra and $C$ a coalgebra. Denote by $J(A)$ (or simply by $J$ when no confusion arises) the Jacobson radical of $A$, by  $(C_n)_{n\geq0}$ the coradical filtration of $C$.  Write $(-)^\perp$ for annihilators under the natural
	pairing between $C$ and its dual algebra $C^*$. Thus $C_n=\bigl(J(C^*)^{n+1}\bigr)^\perp$. 
	
	The notation for the Hopf algebra $H$ is standard: $\Delta$, $\eps$ and $S$ denote the coproduct, the counit and the antipode.
	We use Sweedler notation for the coproduct and coaction. If $H$ is finite-dimensional and $I$ is a Hopf ideal of $H$, then $I^\perp$ is a Hopf subalgebra of $H^*$. Write $G(H)$ for the set of group-like elements of $H$ and $E_H(k)$ for the injective hull of the trivial left $H$-comodule. For $g,h\in G(H)$,  $\Pp_{g,h}(H)=\{x\in H\mid \Delta(x)=x\otimes h+g\otimes x\}$. In particular, $\Pp(H):=\Pp_{1,1}(H)$ is called the set
	of primitive elements. For a subspace $U$ of $H$, denote by $U^r$ the linear span of products of $r$ elements of $U$, by $(U)$ the two-sided ideal generated by $U$, by $\langle U\rangle_{\mathrm{alg}}$ and by $\langle U\rangle_{\mathrm{Hopf}}$ the subalgebra and the Hopf subalgebra generated by $U$, respectively. Let $B(H)$ be
	the subalgebra of $H$ generated by $H_0$.
	If the antipode of $H$ is injective, $B(H)$ is a Hopf
	subalgebra.
	
	A coalgebra is called pointed if all  simple
	subcoalgebras are one dimensional. A finite-dimensional
	Hopf algebra $H$ is basic if all simple modules are
	one dimensional. Over an arbitrary field, we say that $H$ is split basic
	if $H/J(H)\simeq k^s$ for some positive integer $s$. $H$ has the Chevalley property if $H_0$ is a Hopf subalgebra.

	\subsection{Coefficient coalgebras}
	
	Let $C$ be a coalgebra and $M$ a finite-dimensional left $C$-comodule with a fixed basis $m_1,\ldots,m_d$. There exist $c_{i,j}\in C$ such that $\delta(m_j)=\sum_i c_{ji}\otimes m_i$, where 
	\[
	\Delta(c_{ij})=\sum_r c_{ir}\otimes c_{rj},
	\qquad
	\eps(c_{ij})=\delta_{ij}.
	\]
	The subspace $\cf(M)$ of $C$ spanned by the
	coefficients $c_{ji}$ is a subcoalgebra of $C$. $\cf(M)$ is called the coefficient
	coalgebra of $M$. The coefficient coalgebra is independent of the choice
	of basis and depends only on the isomorphism class of $M$.
	
	Suppose that $A$ is a finite-dimensional algebra and $L$ is a left $A$-module with $d=\dim L>0$, equipped with the action map
	$\rho:A\to\End_k(L)$.
	With respect to a basis of $L$, write $f_{ij}(a)=\rho(a)_{ij}$ for $a\in A$. Then 
	$\cf(L)$ in $A^*$ is spanned by coefficients $f_{ij}$. Moreover, 
	\begin{equation}\label{eq:comatrix}
		\Delta(f_{ij})=\sum_r f_{ir}\otimes f_{rj},
		\qquad
		\eps(f_{ij})=\delta_{ij}.
	\end{equation}
	If $\rho$ is surjective, the $d^2$ matrix coefficients
	are linearly independent. The following properties of coefficient coalgebras are standard; 
	see, e.g., \cite[Sections 2.5 and 3.1]{DNR}.
	
	\begin{proposition}\label{prop:coradical-coefficients}
		Let $A$ be a finite-dimensional algebra over $k$, and $\Irr(A)$ a set of representatives of the
		isomorphism classes of simple left $A$-modules. Then
		\[
		(A^*)_0=\bigoplus_{L\in\Irr(A)}\cf(L).
		\]
	\end{proposition}
	
	\begin{lemma}\label{lem:coefficient-coalgebra}
		Let $C$ be a coalgebra and $E$ a finite-dimensional left
		coideal of $C$. Then $\operatorname{cf}(E)$ is the smallest
		subcoalgebra $D$ of $C$ such that $\Delta(E)\subseteq D\otimes E$. 
		In particular, $E\subseteq\operatorname{cf}(E)$.
		Moreover, for a linear map $f:C\to C$, $f|_{\operatorname{cf}(E)}=\id$ if and only if $(f\otimes\operatorname{id})\Delta(e)=\Delta(e)$ for all $e\in E$.
	\end{lemma}

	\subsection{Yetter--Drinfeld modules and Nichols algebras}
	
	Let $\Gamma$ be a finite abelian group and $\YD$ the category
	of left-left Yetter--Drinfeld modules over $k\Gamma$.
	Write $\widehat\Gamma=\Hom(\Gamma,k^\times)$.
	
	Let $V\in\YD$ be finite-dimensional, and suppose that
	the action of $\Gamma$ on $V$ is semisimple. Then
	\[
	V=\bigoplus_{g\in\Gamma,\ \chi\in\widehat\Gamma}V_g^\chi,
	\qquad
	V_g^\chi=
	\{v\in V:\delta(v)=g\otimes v,\
	h\cdot v=\chi(h)v\text{ for all }h\in\Gamma\}.
	\]
	Choose a basis $x_1,\ldots,x_\theta$ of $V$ with $x_i\in V_{g_i}^{\chi_i}$.
	The braiding of $V$ is given by
	\[
	c(x_i\otimes x_j)=q_{ij}x_j\otimes x_i,\quad\text{where}\ q_{ij}=\chi_j(g_i).
	\]
	Thus $(V,c)$ is a braided vector space of diagonal type.
	The collection of Yetter--Drinfeld pairs $(g_i,\chi_i)_{1\leq i\leq\theta}$
	is called a principal realization of $(V,c)$ over $\Gamma$.
	
	\begin{definition}\cite[Definition~2.1]{ASPointed}\label{defiNichols}
		Let $H$ be a Hopf algebra and $V\in{}^{H}_{H}\mathcal{YD}$.
		The Nichols algebra $\B(V)$ over $V$ is a graded
		Hopf algebra $R=\oplus_{n\geq 0}R(n)$ in
		${}^{H}_{H}\mathcal{YD}$ such that
		\begin{align*}
			R(0)=k,\quad R(1)=V,\quad
			R\text{ is generated as an algebra by }R(1),
			\quad \Pp(R)=V.
		\end{align*}
	\end{definition}
	
	The tensor algebra $T(V)=\oplus_{n=0}^{\infty}T^n(V)$
	is a graded Hopf algebra in
	${}^{H}_{H}\mathcal{YD}$, with comultiplication given by
	$\Delta(v)=v\otimes 1+1\otimes v$ for $v\in V$.
	The Nichols algebra $\B(V)$ is isomorphic to
	$T(V)/\mathcal J(V)$, where
	$\mathcal J(V)=\oplus_{n\geq 2}\mathcal J^n(V)$
	is the largest  graded ideal and coideal
	of $T(V)$ in ${}^{H}_{H}\mathcal{YD}$ such that
	$\mathcal J(V)\cap V=0$.
	
	\begin{remark}\label{rmkN-infity}
		The structure of $\B(V)$ depends only on $(V,c)$. Let  $W$ be a vector subspace of $(V, c)$ such that $c(W\otimes W)\subset W\otimes W$. Then $\dim \B(V)=\infty$ if $\dim \B(W)=\infty$. This
		occurs for example, when $V$ contains a non-zero element $v$ such that $c(v\otimes v)=v\otimes v$ in characteristic zero. See \cite{ASPointed} for more details.
	\end{remark}
	\begin{remark}
		For finite-dimensional Nichols algebras of diagonal type
		over an algebraically closed field of characteristic zero,
		a presentation by defining relations is given in
		\cite{AngionoRelations}.
	\end{remark}

	\subsection{Bosonization and liftings}
	Let $R$ be a Hopf algebra in ${}^{H}_{H}\mathcal{YD}$. Write $\Delta_R(r)=r^{(1)}\otimes r^{(2)}$ for the braided coproduct to avoid confusions. The bosonization $R\# H$ is defined as follows: the underlying space $R\# H=R\otimes H$, the product and coproduct are given by  
	\begin{align}
		(r\# g)(s\# h)=r(g_{(1)}\cdot s)\# g_{(2)}h,\quad
		\Delta(r\# g)=r^{(1)}\#(r^{(2)})_{(-1)}g_{(1)}\otimes (r^{(2)})_{(0)}\# g_{(2)}.\label{eqSmash}
	\end{align}
	The map $j:H\rightarrow R\# H, h\mapsto 1\# h,\ \forall h\in H$ is injective and  
	$\pi:R\# H\rightarrow H,r\# h\mapsto \epsilon_R(r)h,\ \forall r\in R, h\in H$ is surjective  satisfying $\pi\circ j=\id_H$.
	Moreover, $R=(R\# H)^{coH}=\{x\in R\# H\mid (\text{id}\otimes\pi)\Delta(x)=x\otimes 1\}$.

	Conversely, if $A$ is a Hopf algebra and $\pi:A\rightarrow H$ is a bialgebra map admitting a bialgebra
	section $j:H\rightarrow A$ such that $\pi\circ j=\id_H$, then $A\simeq R\# H$, where $R=A^{coH}$ is a Hopf algebra in ${}^{H}_{H}\mathcal{YD}$ with 
	\begin{gather}\label{eQ-sh}
		h\cdot r=j(h_{(1)})rj(S_H(h_{(2)})),\quad
		\delta(r)=(\pi\otimes id)\Delta_A(r),\quad
		\Delta_R(r)=r_{(1)}(j S_H(\pi(r_{(2)})))\otimes r_{(3)}.
		%\epsilon_R=\epsilon_A|_R,\;
		%S_R(r)=(j\pi(r_{(1)}))S_A(r_{(2)}).
	\end{gather}
	
	Let $A$ be a finite-dimensional pointed Hopf algebra
	with $G(A)=\Gamma$ and $\pi:\gr A\longrightarrow k\Gamma$ the canonical
	Hopf algebra projection. The diagram of
	$A$ is the graded braided Hopf algebra
	\[
	R=(\gr A)^{\operatorname{co}\pi}
	=
	\{a\in\gr A:
	(\Id\otimes\pi)\Delta(a)=a\otimes1\}\ \text{in }\ \YD,
	\]
	so that $\gr A\simeq R\#k\Gamma$.
	Moreover, $R$ is strictly graded, namely, $R(0)=k1$ and  $R(1)=\Pp(R)$.
	The space $R(1)$ is called the infinitesimal
	braiding of $A$.

	For a given $V\in\YD$, $A$ is a lifting of 
	$\B(V)\#k\Gamma$ if $\gr A\simeq\B(V)\#k\Gamma$. For $g\in\Gamma$ and $\chi\in\widehat\Gamma$, put
	\[
	\Pp_{g,1}^{\chi}(A)
	=
	\{x\in\Pp_{g,1}(A):
	hxh^{-1}=\chi(h)x
	\text{ for all }h\in\Gamma\}.
	\]
	
	For any $x_i\in V_{g_i}^{\chi_i}$,  by
	\cite[Lemma~5.1]{ASPointed}
	one can choose a skew-primitive representative
	$\widetilde{x}_i\in A_1$ such that
	\[
	h\widetilde{x}_ih^{-1}
	=\chi_i(h)\widetilde{x}_i,\qquad
	\Delta(\widetilde{x}_i)
	=\widetilde{x}_i\otimes1+g_i\otimes\widetilde{x}_i
	\quad(h\in\Gamma),
	\]
	and whose image in $A_1/A_0\simeq V\#k\Gamma$
	is $x_i\#1$. Write
	$\iota:V\to A_1$ for the linear map
	$\iota(x_i)=\widetilde{x}_i$ and also $x_i$ for $\iota(x_i)$ in $A$ when no confusion arises. We retain $\iota$
	when a subspace of $V$ and its image in $A$ occur in
	the same calculation.
	
	For $g\in\Gamma$ and $\chi\in\widehat\Gamma$, set
	\[
	E_{g,\chi}
	=
	\operatorname{span}\{x_i:g_i=g,\ \chi_i=\chi\}
	\subseteq A_1.
	\]
	The natural map $A_1\to A_1/A_0$ identifies the image
	of $E_{g,\chi}$ with $V_g^\chi\#1$.
	
	Assume that $k$ is algebraically closed of characteristic
	zero and that $\gr A\simeq\B(V)\#k\Gamma$ is
	finite-dimensional. Then by \cite[Lemma~5.1]{ASPointed},
	\[
	\Pp_{g,1}(A)=k(1-g)\oplus\bigoplus_{\chi\ne\eps}E_{g,\chi},
	\]
	where
	\[
	E_{g,\chi}=\iota(V_g^\chi)=\Pp_{g,1}^{\chi}(A)
	\quad(\chi\ne\eps),
	\qquad
	\Pp_{g,1}^{\eps}(A)=k(1-g).
	\]

	\section{The Hopf coradical}\label{sec:structure}
	We characterize the Hopf coradical and relate it to semisimple Hopf subalgebras and injective comodules.
	
	\subsection{A criterion for generation by the coradical}We begin with the following result concerning a subspace that generates an algebra modulo a square-zero ideal.
	\label{sec:basic}
	
\begin{lemma}\label{lem:squarezero}
		Let $R$ be a unital algebra over a field $k$,  $M$ a
		two-sided ideal with $M^2=0$, and $C$ a linear subspace
		of $R$ containing $1$. If the canonical map $\pi: C\to R/M$ is surjective, then $\langle C\rangle_{\mathrm{alg}}=C^3$.
	\end{lemma}
\begin{proof}
	Let $Q=R/M$ and $D=\langle C\rangle_{\mathrm{alg}}$.
	Since $\pi:C\to Q$ is surjective, there exists a linear
	section $\sigma:Q\to C$ such that
	$\sigma(1)=1$ and $\pi\sigma=\id_Q$.
	$M$ admits a $Q$-bimodule structure by setting
	\[
	q\cdot m=\sigma(q)m,\qquad
	m\cdot q=m\sigma(q),
	\qquad\text{for }q\in Q,\ m\in M.
	\]
	The $Q$-actions are well defined; indeed, if $r\in R$
	represents $q$, then $r-\sigma(q)\in M$, and since
	$M^2=0$, $rm=\sigma(q)m$ and $mr=m\sigma(q)$.
	
	Set $C_M=C\cap M$ and define
	\[
	\alpha:Q\times Q\longrightarrow M,\qquad
	\alpha(a,b)=\sigma(a)\sigma(b)-\sigma(ab).
	\]
	The map $\alpha$ is well-defined because the canonical
	quotient map $R\to Q$ sends both $\sigma(a)\sigma(b)$
	and $\sigma(ab)$ to $ab$.
	
	Let $N$ be the $Q$-subbimodule of $M$ generated by
	$C_M$ and the elements $\alpha(a,b)$ for all $a,b\in Q$.
	By construction of the  $Q$-bimodule on $M$,
	$N$ is stable under the $R$-action on $M$; hence $N$ is
	a two-sided ideal of $R$. Moreover, $N\subseteq D$;
	indeed, $C_M\subseteq C\subseteq D$,
	$\alpha(a,b)\in C^2+C\subseteq D$, and the $Q$-actions
	on $N$ are induced by multiplication by
	$\sigma(Q)\subseteq C\subseteq D$.

	Let $\overline R=R/N$ and let
	$\rho:R\to\overline R$ be the canonical quotient map.
	Define $\overline\sigma=\rho\circ\sigma:Q\to\overline R$.
	Since $\alpha(a,b)\in N$ for $a,b\in Q$, we have $\overline\sigma(a)\overline\sigma(b)=\overline\sigma(ab)$.
	That is, $\overline\sigma$ is a unital algebra homomorphism,
	and so $\overline\sigma(Q)$ is a subalgebra of $\overline R$.
	For every $c\in C$,
	\[
	c-\sigma(\pi(c))\in C\cap M=C_M\subseteq N.
	\]
	Then the image of $C$ in $\overline R$
	is precisely $\overline\sigma(Q)$ and so the image of
	$D=\langle C\rangle_{\mathrm{alg}}$  is also
	$\overline\sigma(Q)$. As $N\subseteq D$, one has $D=\sigma(Q)+N$.

	It remains to show that $N\subseteq C^3$. Since $1\in C$,
	we have $C\subseteq C^2$, and hence $\alpha(a,b)\in C^2$
	for all $a,b\in Q$. For $c_M\in C_M$ and $a,b\in Q$,
	$\sigma(a)c_M\sigma(b)\in C^3$. For the other generators,
	associativity and $M^2=0$ yield, for all $x,a,b,y\in Q$,
	\[
	\sigma(x)\alpha(a,b)\sigma(y)
	=\alpha(xa,b)\sigma(y)-\alpha(x,ab)\sigma(y)
	+\alpha(x,a)\sigma(by).
	\]
	Since $\alpha(u,v)\in C^2$ and $\sigma(Q)\subseteq C$,
	every term on the right lies in $C^2\cdot C=C^3$,
	which implies that $N\subseteq C^3$. As
	$\sigma(Q)\subseteq C\subseteq C^3$, we obtain
	$D=\sigma(Q)+N\subseteq C^3$, and the reverse inclusion
	follows from $C\subseteq D$. Hence $D=C^3$.
\end{proof}

	\begin{theorem}\label{thm:nilpotent-generation}
		Let $R$ be a unital algebra, $I$ a nilpotent two-sided
		ideal of $R$ and $\pi:R\to R/I$ be the quotient map. 
		Suppose that $C$ is a subspace of $R$ containing $1$ such
		that the restriction $\psi=\pi|_C:C\to R/I$ is surjective.
		Then $R=\langle C\rangle_{\mathrm{alg}}$ if and only if $R=C^3+I^2$.
		If these conditions hold and $L=C^3\cap I$, then, for
		any $\ell\geq 1$ such that $I^\ell=0$,
		\begin{equation}\label{eq:radical-filtration}
			I^r=L^r+I^{r+1}
			=\sum_{j=r}^{\ell-1}L^j,\qquad 1\leq r<\ell.
		\end{equation}
	\end{theorem}
	
	\begin{proof}
		Suppose that $R=\langle C\rangle_{\mathrm{alg}}$.
		Let $\overline R=R/I^2$, $\overline I=I/I^2$ and
		$\overline C=(C+I^2)/I^2$.
		Then $\overline I^2=0$, and the map
		$\overline\psi:\overline C\to R/I$,
		$\overline\psi(c+I^2)=\psi(c)$, is well-defined and
		surjective. Moreover, $\overline C$ contains $1$ and
		$\overline R=\langle\overline C\rangle_{\mathrm{alg}}$.
		By Lemma~\ref{lem:squarezero},
		$\overline R=\overline C^3$ and so $R=C^3+I^2$.
		
		Conversely, suppose that $R=C^3+I^2$, and set
		$L=C^3\cap I$. For any $z\in I\subset R$, there are $c\in C^3$ and $u\in I^2$ such that $z=c+u$. Since $I^2\subseteq I$,
		one has $c=z-u\in C^3\cap I=L$. Thus $I=L+I^2$.
		
		For $1\leq r<\ell$, expanding $I^r=(L+I^2)^r$,
		every term other than $L^r$ belongs to $I^{r+1}$, which means that $I^r=L^r+I^{r+1}$. Iterating $I^r=L^r+I^{r+1}$ and using
		$I^\ell=0$, we obtain \eqref{eq:radical-filtration}.
		Since $L\subseteq C^3$, $I\subseteq\langle C\rangle_{\mathrm{alg}}$.
		From the assumption that $\psi:C\to R/I$ is surjective, $R=\langle C\rangle_{\mathrm{alg}}$ follows.
	\end{proof}

	\begin{example}\label{ex:length-three}
		In general, the exponent $3$ in Lemma~\ref{lem:squarezero} cannot be replaced
		by $2$. Take $Q=ke_1\oplus ke_2$, where
		$e_ae_b=\delta_{ab}e_a$. Let $M$ be the $Q$-bimodule with
		basis $m_{ij}$, $1\le i,j\le2$ with the action given by 
		$e_am_{ij}=\delta_{ai}m_{ij}$ and
		$m_{ij}e_b=\delta_{jb}m_{ij}$.
		
		Define an algebra $R=Q\oplus M$ where the multiplication is given as follows: $(a,m)(b,n)=(ab,an+mb)$. Then $M^2=0$.
		Let $m=\sum_{i,j}m_{ij}$ and
		$C=\operatorname{span}\{e_1,e_2,m\}$.
		
		The subspace $C^2\cap M$ is spanned by the two row sums
		$e_i m$ and the two column sums $m e_j$. The four elements $e_1m,e_2m,me_1,me_2$ satisfy
		the unique linear relation $e_1m+e_2m=me_1+me_2$. Thus $\dim C^2=5$. On the other hand,
		$e_i m e_j=m_{ij}$, so $C^3=R$.
		
		Hence the exponent $3$ cannot be replaced by $2$ for
		arbitrary subspaces $C$. However, this example does not settle
		whether $C^2$ suffices when $C$ is the coradical of a
		finite-dimensional Hopf algebra.
	\end{example}

	\begin{lemma}\label{lem:surject-equivalence}
		Let $H$ be a finite-dimensional Hopf algebra and set $A=H^*$. Then the canonical map $H_0\longrightarrow H/J(H)$ is surjective precisely when $J(A)\cap A_0=0$.
	\end{lemma}
	\begin{proof}
		Let $\pi:H\to H/J(H)$ be the quotient map, and set
		$T=\pi|_{H_0}:H_0\to H/J(H)$.
		Since $(H/J(H))^*\simeq J(H)^\perp=A_0$, the transpose $T^*:(H/J(H))^*\longrightarrow H_0^*$ is exactly the restriction map
		\[
		\pi^*:A_0\longrightarrow H_0^*,
		\qquad f\longmapsto f|_{H_0}.
		\]
		Indeed, if $f\in A_0$ corresponds to
		$\bar f\in(H/J(H))^*$, then $T^*(\bar f)(h)=\bar f(\pi(h))=f(h)$ for $h\in H_0$. 
		Hence,
		\[
		\ker\pi^*=\{f\in A_0:f|_{H_0}=0\}=A_0\cap H_0^\perp=A_0\cap J(A).
		\]
		Consequently, $T$ is surjective if and only if $\ker\pi^*=0$, or equivalently, $A_0\cap J(A)=0$.
	\end{proof}

	\begin{lemma}\label{lem:surject}
		Let $H$ be a finite-dimensional Hopf algebra over $k$. If $(H^*)_0$ is a semisimple subalgebra of $H^*$, then $H_0\to H/J(H)$ is surjective. In particular, the hypothesis is satisfied if $H$ is split basic and $\operatorname{char}k\nmid |G(H^*)|$, or if $\operatorname{char}k=0$ and $H^*$ has the Chevalley property.
	\end{lemma}
	
	\begin{proof}
		Let $A=H^*$. The ideal $J(A)\cap A_0$ is a nilpotent two-sided ideal of $A_0$, hence zero, since $A_0$ is semisimple. By Lemma \ref{lem:surject-equivalence}, $H_0\to H/J(H)$ is surjective.
		
		When $H$ is split basic, $A$ is pointed and $A_0=kG(A)$, which is semisimple by Maschke's theorem. If $A$ has the Chevalley property in characteristic zero, $A_0$ is a cosemisimple Hopf subalgebra of $A$, hence semisimple by \cite[Theorem 16.1.2]{RadfordBook}.
	\end{proof}
	
	\begin{theorem}\label{thm:four}
		Let $H$ be a finite-dimensional Hopf algebra over a field
		$k$. Suppose that $H_0\to H/J(H)$ is surjective. Then the following are equivalent:
		\begin{itemize}
			\item[(i)] $H$ is generated by $H_0$, that is, $H=B(H)$;
			\item[(ii)] $H=H_0^3+J(H)^2$;
			\item[(iii)] $H=H_0^4+J(H)^2$.
		\end{itemize}
		In particular, the equivalences hold if $H$ is split basic
		and $\operatorname{char}k\nmid |G(H^*)|$, or if
		$\operatorname{char}k=0$ and $H^*$ has the Chevalley property.
	\end{theorem}
	
	\begin{proof}
		Set $J=J(H)$. Since $J$ is nilpotent and $1\in H_0$,
		the equivalence of $(i)$ and $(ii)$ follows by Theorem~\ref{thm:nilpotent-generation}.
		Since $H_0^3\subseteq H_0^4\subseteq B(H)$, $(i)$ implies
		$(iii)$ by $(ii)$.
		Conversely, $(iii)$ gives $H=B(H)+J^2$, so the map
		$B(H)\to H/J$ is surjective. As $B(H)^3=B(H)$, by Theorem~\ref{thm:nilpotent-generation}, we obtain $H=B(H)$.
		The remaining assertions follow from
		Lemma~\ref{lem:surject}.
	\end{proof}

	\begin{corollary}\label{cor:length}
		Under the hypotheses of Theorem~\ref{thm:four}, suppose that
		$B(H)=H$. Let $J=J(H)$ and $L=H_0^3\cap J$.
		If $J^\ell=0$ for some $\ell\geq2$, then
		\[
		J^r=\sum_{j=r}^{\ell-1}L^j\quad(1\leq r<\ell),
		\qquad
		H=H_0^{3(\ell-1)}.
		\]
		If $J=0$, then $H=H_0$.
	\end{corollary}
	
	\begin{proof}
		If $J=0$, then clearly $H=H_0$. Otherwise the formula for $J^r$ follows from
		\eqref{eq:radical-filtration}. Since $H_0$ maps onto
		$H/J$, we have $H=H_0+J$. Since $1\in H_0$, one has $L^j\subseteq H_0^{3j}\subseteq H_0^{3(\ell-1)}$
		for $1\leq j<\ell$.
		As a result, $H=H_0^{3(\ell-1)}$.
	\end{proof}
	
	\begin{corollary}\label{cor:rank}
		Under the hypotheses of Theorem~\ref{thm:four}, the
		multiplication map
		\[
		m_3:H_0^{\otimes3}\longrightarrow H/J(H)^2,
		\qquad
		a\otimes b\otimes c\longmapsto abc+J(H)^2
		\]
		has rank $\dim H/J(H)^2$ if and only if $B(H)=H$.
	\end{corollary}
	
	\begin{proof}
		Observe that the image of $m_3$ is
		$(H_0^3+J(H)^2)/J(H)^2$. Then $m_3$ is surjective
		if and only if $H=H_0^3+J(H)^2$, which is equivalent
		to $B(H)=H$ by Theorem~\ref{thm:four}.
	\end{proof}
	
	(CG) can be also characterized in terms of the injective hull of the
	trivial comodule.
	
	\begin{theorem}\label{thm:injective-product}
		Let $H$ be a finite-dimensional Hopf algebra over
		any field, $E_H(k)$ the injective hull
		of the trivial left $H$-comodule and $\C(H)$
		the Hopf subalgebra of $H$ generated by $\operatorname{cf}(E_H(k))$. Then
		\[
		H=H_0\C(H)=\C(H)H_0,
		\]
		Furthermore, $H$ is generated by $H_0$ exactly when $\C(H)\subseteq B(H)$.
	\end{theorem}
	
	\begin{proof}
		Since $H$ is injective as a left $H$-comodule,
		the inclusion $k1\subseteq H$ can extend to a left
		$H$-comodule map
		$\psi:E_H(k)\to H$. As $k$ is essential in
		$E_H(k)$, $\psi$ is injective. Put
		$E=\operatorname{im}\psi$. Then $E$ is a left
		coideal of $H$, and
		$\operatorname{cf}(E)=\operatorname{cf}(E_H(k))$.
		
		By Lemma~\ref{lem:coefficient-coalgebra},
		$E\subseteq\operatorname{cf}(E)$. If
		$K$ is a Hopf subalgebra containing $E$, then
		$\Delta(E)\subseteq K\otimes E$, and hence
		$\operatorname{cf}(E)\subseteq K$. As a result,
		\[
		\C(H)=\langle E\rangle_{\mathrm{Hopf}}
		=\langle\operatorname{cf}(E)\rangle_{\mathrm{Hopf}}.
		\]
		
		Since $H$ is finite-dimensional, it is
		co-Frobenius. By~\cite[Proposition 10.9.10(d)]{RadfordBook},
		$H=H_0E$. Thus, 
		\[
		H=H_0E\subseteq H_0\C(H)\subseteq H,
		\]
		so $H=H_0\C(H)$. Since $S(H_0)=H_0$ and
		$S(\C(H))=\C(H)$, one has
		\[
		H=S(H)=S(\C(H))S(H_0)=\C(H)H_0.
		\]
		
		If $\C(H)\subseteq B(H)$, then
		$H=H_0\C(H)\subseteq B(H)$, and hence $H=B(H)$.
		The converse follows from $\C(H)\subseteq H$.
	\end{proof}
	\subsection{The Jacobson radical of the dual Hopf algebra and Skew-primitive elements}
	\label{subsec:dual-radical}
	Let $H$ be a finite-dimensional Hopf algebra and $\Imax(H)$ the largest Hopf ideal of $H$ contained in $J(H)$.
	The following result is due to Chen and Hiss.
	
	\begin{lemma}[{\cite[Corollary 2.3]{CH}}]\label{lem:largest-coideal-in-radical}
		The ideal $\Imax(H)$ is also the largest coideal contained in $J(H)$.
	\end{lemma}

	\begin{lemma}\label{lem:dual-radical-annihilator}
		If $H$ is finite-dimensional, then $\Imax(H^*)=B(H)^\perp$.
	\end{lemma}
	
	\begin{proof}
		By definition, $H_0\subseteq B(H)$ and so $B(H)^\perp\subseteq H_0^\perp=J(H^*)$.
		 $B(H)^\perp$ is a Hopf ideal of $H^*$ contained in $J(H^*)$.
		Assume that $I$ is a Hopf ideal contained in $J(H^*)$ such that $B(H)^\perp\subseteq I$.
		Then $I^\perp$ is a Hopf subalgebra of $H$ containing
		$H_0$. Hence $B(H)\subseteq I^\perp$, and
		$I\subseteq B(H)^\perp$. Consequently, $\Imax(H^*)=B(H)^\perp$.
	\end{proof}

	For an $H^*$-module $M$, the Hopf image of $M$ is the quotient of
	$H^*$ by the largest Hopf ideal contained in
	$\operatorname{Ann}_{H^*}(M)$ \cite{BB}.
	The module is inner faithful if this Hopf ideal is zero.
	
	Let $M_{\mathrm{sim}}$ be the direct sum of representatives of
	the isomorphism classes of simple $H^*$-modules.
	Since the Jacobson radical is the intersection of the
	annihilators of simple modules, we have
	\[
	\operatorname{Ann}_{H^*}(M_{\mathrm{sim}})=J(H^*).
	\]
	By Lemma~\ref{lem:dual-radical-annihilator},
	one can describe the Hopf image as follows:
	
	\begin{proposition}
		\label{prop:hopf-image-simple-modules}
		The Hopf image of $M_{\mathrm{sim}}$ is canonically
		isomorphic to $B(H)^*$. In particular,
		$H=B(H)$ if and only if $M_{\mathrm{sim}}$ is
		inner faithful.
	\end{proposition}
	
	\begin{corollary}
		\label{cor:dual-chevalley-radical}
		If $H$ is finite-dimensional, then $H$ has the Chevalley property
		if and only if $\Imax(H^*)=J(H^*)$.
	\end{corollary}
	
	\begin{proof}
		Note that $B(H)=H_0$ if and only if $H_0$ is a Hopf subalgebra. The assertion follows from Lemma \ref{lem:dual-radical-annihilator} and $J(H^*)=H_0^\perp$.
	\end{proof}

	\begin{lemma}\label{lem:pointed-detect}
		Let $A$ be a finite-dimensional pointed Hopf algebra over a field
		$k$ with $\operatorname{char}k\nmid |G(A)|$. Then $A^*$ is generated by the coradical if and only if $J(A)\cap\Pp_{g,1}(A)=0$ for every $g\in G(A)$.
	\end{lemma}

	\begin{proof}
		By Lemma~\ref{lem:dual-radical-annihilator}, $B(A^*)=A^*$ if and only if $\Imax(A)=0$.
		If $z\in J(A)\cap\Pp_{g,1}(A)$, then clearly $(z)$ is a Hopf
		ideal contained in $J(A)$ and so $(z)\subseteq\Imax(A)$. Consequently,
		$\Imax(A)=0$ implies $z=0$.
		
		Conversely, let $0\ne I\subseteq J(A)$ be a Hopf ideal,
		and let $\pi:A\to A/I$ be the quotient map. By Maschke's
		theorem, $A_0=kG(A)$ is semisimple. Since $I\cap A_0$
		is a nilpotent two-sided ideal of $A_0$, it follows that
		$I\cap A_0=0$. In particular, $\pi\mid_{A_0}$ is injective.
		
		By \cite[Theorem 5.3.1]{Montgomery}, $\pi\mid_{A_1}$ cannot be
		injective; otherwise $\pi$ would be injective.  Observe that $A/I$ is pointed
		and $(A/I)_0=\pi(A_0)$. The induced
		map
		\[
		\overline\pi:A_1/A_0\longrightarrow
		(A/I)_1/(A/I)_0
		\]
		has nonzero kernel.  
		
		For $a,b\in G(A)$,  denote by $\overline{\Pp}_{a,b}$ the image of $\Pp_{a,b}(A)$ in $A_1/A_0$.
		By Taft--Wilson Theorem, 
		\[
		A_1/A_0
		=
		\bigoplus_{a,b\in G(A)}
		\overline{\Pp}_{a,b}.
		\]
		The induced map $\overline\pi$ sends
		$\overline{\Pp}_{a,b}$ into the component indexed by
		$(\pi(a),\pi(b))$ in the corresponding decomposition
		of $(A/I)_1/(A/I)_0$.
		Since $\pi\mid_{A_0}$ is injective, distinct pairs
		$(a,b)$ give distinct pairs $(\pi(a),\pi(b))$.
		There exist $a,b\in G(A)$ and
		$z\in\Pp_{a,b}(A)\setminus A_0$ such that $\pi(z)\in\pi(A_0)$.

		Since $\Pp_{\pi(a),\pi(b)}(\pi(A_0))=k\bigl(\pi(a)-\pi(b)\bigr)$,  there is $\lambda\in k$ such that $\pi(z)=\lambda\bigl(\pi(a)-\pi(b)\bigr)$. Since $z\notin A_0$,
		\[
		0\ne z-\lambda(a-b)\in I\cap\Pp_{a,b}(A),
		\]
		which implies that $0\ne \bigl(z-\lambda(a-b)\bigr)b^{-1} \in I\cap\Pp_{ab^{-1},1}(A)$.
	\end{proof}

	\begin{corollary}\label{cor:equal-type}
		Let $A$ be a finite-dimensional pointed Hopf algebra with abelian coradical $k\Gamma$ over an algebraically closed field of characteristic zero. Then $A^*$ is generated by the coradical if and only if $J(A)\cap\Pp_{g,1}^{\chi}(A)=0$
		for all $g\in\Gamma$ and $\chi\in\widehat\Gamma\setminus\{\varepsilon\}$.
	\end{corollary}
	\begin{proof}
		Since $\Gamma$ is finite abelian and $\operatorname{char}k=0$, 
		\[
		\Pp_{g,1}(A)
		=k(1-g)\oplus\bigoplus_{\chi\ne\varepsilon}
		\Pp_{g,1}^{\chi}(A).
		\]
		Moreover, $J(A)\cap k(1-g)=0$, since $k(1-g)\subseteq A_0$
		and $A_0=k\Gamma$ is semisimple. For $x\in J(A)\cap\Pp_{g,1}(A)$,
		there exist $\lambda\in k$ and
		$x_\chi\in\Pp_{g,1}^{\chi}(A)$ such that
		\[
		x=\lambda(1-g)+\sum_{\chi\ne\varepsilon}x_\chi.
		\]
		For every $\chi\in\widehat\Gamma$, we have
		\[
		j_1:=\frac{1}{|\Gamma|}\sum_{h\in\Gamma}
		\chi(h)^{-1}hxh^{-1}\in J(A),
		\]
		since $J(A)$ is a two-sided ideal. If $\chi\ne\varepsilon$,
		then $j_1=x_\chi$; otherwise $j_1=\lambda(1-g)$.
		Since $J(A)\cap k(1-g)=0$, we obtain
		$x_\chi\in J(A)\cap\Pp_{g,1}^{\chi}(A)$ for all
		$\chi\ne\varepsilon$. Consequently,
		\[
		J(A)\cap\Pp_{g,1}(A)
		=\bigoplus_{\chi\ne\varepsilon}
		\bigl(J(A)\cap\Pp_{g,1}^{\chi}(A)\bigr).
		\]
		The assertion follows from Lemma~\ref{lem:pointed-detect}.
	\end{proof}
	
	Let $H$ be a finite-dimensional basic Hopf algebra over an
	algebraically closed field of characteristic zero. For $g\in G(H^*)$, set  
	\[
	\mathcal D_g:=\{d\in H^*:d(ab)=d(a)\varepsilon(b)+g(a)d(b)\text{ for all }a,b\in H\}.
	\]
	Then $\mathcal D_g=\Pp_{g,1}(H^*)$.
	
	\begin{corollary}\label{cor:rank-criteria}
		Let $H$ be a finite-dimensional basic Hopf algebra over an
		algebraically closed field of characteristic zero.
		For each $g\in G(H^*)$, let
		$\{d_{g,1},\ldots,d_{g,r_g}\}$ and $\{c_1,\ldots,c_s\}$ be bases of $\mathcal D_g$ and $H_0$, respectively.
		Set $M_g=(d_{g,j}(c_i))_{i,j}$. Then the following are
		equivalent:
		\begin{enumerate}
			\item $B(H)=H$;
			\item $\rank M_g=r_g$ for every $g\in G(H^*)$;
			\item $\rank m_3=\dim H/J(H)^2$, where $m_3$ is defined in Corollary \ref{cor:rank}.
		\end{enumerate}
	\end{corollary}
	
	\begin{proof}
		The restriction map $\mathcal D_g\to H_0^*$ has matrix
		$M_g$ with respect to the chosen bases, and its kernel is
		\[
		\mathcal D_g\cap H_0^\perp
		=\mathcal D_g\cap J(H^*).
		\]
		Since $H^*$ is pointed and $\operatorname{char}k=0$,
		Lemma~\ref{lem:pointed-detect} shows that (1) and (2)
		are equivalent.
		
		Recall that the image of $m_3$ is
		$(H_0^3+J(H)^2)/J(H)^2$. Then (3) is equivalent to
		$H=H_0^3+J(H)^2$ and hence is equivalent to (1) by Theorem~\ref{thm:four}.
	\end{proof}

	\subsection{Semisimple Hopf subalgebras and the antipode}
	Let $H$ be a finite-dimensional Hopf algebra in characteristic zero. Larson and Radford
	(see \cite[Theorem 16.1.2]{RadfordBook}) show that $H$ is
	semisimple if and only if it is cosemisimple, or equivalently,
	if $S^2=\operatorname{id}_H$.
	\begin{proposition}\label{prop:max-semisimple}
		Let $H$ be a finite-dimensional Hopf algebra over a field of
		characteristic zero and set
		\[
		H_{\mathrm{ss}}
		=\{h\in H:(\Id\otimes S^2)\Delta(h)=\Delta(h)\}.
		\]
		Then $H_{\mathrm{ss}}$ is the largest semisimple Hopf subalgebra of $H$ contained in $H_0$.
	\end{proposition}
	\begin{proof}
		Let $F:=S_{H^*}^2$, and  $L$ the ideal of $H^*$ generated by $F(b)-b$ for all $b\in H^*$.
		Then $F$ is a Hopf algebra automorphism of $H^*$,  and 
		\[
		\Delta(F(b)-b)
		=F(b_1)\otimes(F(b_2)-b_2)
		+(F(b_1)-b_1)\otimes b_2,
		\]
		and $\varepsilon(F(b)-b)=0$.
		Moreover, $S_{H^*}$ commutes with $F$, and so
		\[
		S_{H^*}(F(b)-b)
		=F(S_{H^*}(b))-S_{H^*}(b)\in L.
		\]
		As a result, $L$ is a Hopf ideal. Observe that
		\[
		(F(b)-b)c
		=(F(bc)-bc)-F(b)(F(c)-c).
		\]
		The left ideal generated by $F(b)-b$ for all $b\in H^*$
		is also a right ideal. In particular,
		\[
		L=\operatorname{span}\{a(F(b)-b):a,b\in H^*\}.
		\]
		
		Under the canonical pairing between $H^*$ and $H$, we have
		\[
		\begin{aligned}
			\langle a(F(b)-b),h\rangle
			&=\langle a\otimes(F(b)-b),\Delta(h)\rangle\\
			&=\langle a\otimes b,
			(\operatorname{id}_H\otimes S_H^2)\Delta(h)-\Delta(h)
			\rangle.
		\end{aligned}
		\]
		Since the pairing between $H^*\otimes H^*$ and $H\otimes H$
		is nondegenerate, it follows that $h\in L^\perp$ if and only if  $(\operatorname{id}_H\otimes S_H^2)\Delta(h)=\Delta(h)$.
		Thus $H_{\mathrm{ss}}=L^\perp$ is a Hopf subalgebra of $H$.
		
		For $h\in H_{\mathrm{ss}}$,  
		\[
		S_H^2(h)=\varepsilon(h_1)S_H^2(h_2)
		=\varepsilon(h_1)h_2=h.
		\]
		Thus $H_{\mathrm{ss}}$ is semisimple and
		cosemisimple and so $H_{\mathrm{ss}}\subseteq H_0$.
		
		Finally, let $T$ be a semisimple Hopf subalgebra of $H$.
		Then $S_T^2=\operatorname{id}_T$ and furthermore $S_H^2|_T=\operatorname{id}_T$. As $\Delta(T)\subseteq
		T\otimes T$, 
		\[
		(\operatorname{id}_H\otimes S_H^2)\Delta(t)=\Delta(t), \qquad t\in T.
		\]
		Consequently, $T\subseteq H_{\mathrm{ss}}$.
	\end{proof}
	
	\begin{proposition}\label{prop:coideal-antipode}
		Let $H$ be a finite-dimensional Hopf algebra over a field of
		characteristic zero and $E$ a left coideal of $H$. Suppose that $K$ is the Hopf subalgebra of $H$ generated by $E$.
		Then $K$ is semisimple if and only if $S^2(e)=e$ for every $e\in E$.
	\end{proposition}
	\begin{proof}
		If $K$ is semisimple, then $S_K^2=\operatorname{id}_K$, and hence $S_H^2(e)=e$ for all $e\in E$.
		
		Conversely, suppose that $S_H^2(e)=e$ for all $e\in E$.
		Then from $\Delta(E)\subseteq H\otimes E$, for any $e\in E$,  
		\[
		\Delta(e)
		=\Delta(S_H^2(e))
		=(S_H^2\otimes S_H^2)\Delta(e)
		=(S_H^2\otimes\operatorname{id}_E)\Delta(e).
		\]
		By definition,  $S_H^2\mid_{\operatorname{cf}(E)}=\id$.
		
		By Lemma \ref{lem:coefficient-coalgebra}, $E\subseteq\operatorname{cf}(E)\subseteq K$ and hence $K$ is generated by  $\operatorname{cf}(E)$. Since $S_H^2$ commutes with $S_H$ and is an algebra
		automorphism, it follows that $S_H^2|_K=\operatorname{id}_K$ and hence $K$ is semisimple.
	\end{proof}
	
	\begin{remark}\label{rem:semisimple-normal}
		For a finite-dimensional Hopf algebra, the largest semisimple Hopf subalgebra need not be normal. Let $H_4$ be Sweedler's Hopf algebra generated by $g,x$ over a field of characteristic zero, subject to the relations $g^2=1$, $x^2=0$, $gx=-xg$, and
		$\Delta(x)=x\otimes1+g\otimes x$. The largest semisimple Hopf subalgebra of $H_4$ is $k\langle g\rangle$, while $k\langle g\rangle$ is not normal in $H_4$.
	\end{remark}

	\begin{corollary}\label{cor:chevalley-ss}
		Let $H$ be a finite-dimensional Hopf algebra in characteristic zero. Then $H$ has the Chevalley property if and only if \(H_{\mathrm{ss}}=H_0\), or
		equivalently \(S^2|_{H_0}=\Id\).
	\end{corollary}
	\begin{proof}
		If \(H_0\) is a Hopf subalgebra, \(H_0\) is cosemisimple and hence
		semisimple. By Proposition~\ref{prop:max-semisimple}, $H_{\mathrm{ss}}=H_0$. Conversely, $H_0$ must be a Hopf
		subalgebra and hence $H$ has the Chevalley property. If \(S^2\mid_{H_0}=\id\), then $\Delta(H_0)\subseteq H_0\otimes H_0$ implies
		$H_0\subseteq H_{\mathrm{ss}}$. Consequently, $H_0=H_{\mathrm{ss}}$ follows
		by Proposition~\ref{prop:max-semisimple}.
	\end{proof}

	We end up by determining the coradical of the product of a semisimple
	normal Hopf subalgebra and a Hopf subalgebra with the Chevalley
	property.
	
	\begin{proposition}\label{prop:semisimple-product}
		Let $H$ be a finite-dimensional Hopf algebra over a field of
		characteristic zero. Suppose that $U$ is a semisimple normal Hopf
		subalgebra of $H$, and $V$ is a Hopf subalgebra of $H$ with the
		Chevalley property.
		Then $UV$ is a Hopf subalgebra of $H$ with the Chevalley
		property, and
		\[
		(UV)_0=UV_0,\qquad
		\dim UV=\frac{\dim U\dim V}{\dim(U\cap V)}.
		\]
	\end{proposition}
	\begin{proof}
		Since $U$ is normal, $UV=VU$ is a Hopf subalgebra of $H$.
		Let $U^+=\ker(\varepsilon|_U)$ and 
		$Q=UV/(UV)U^+$. Denote the quotient map by $q:UV\to Q$.
		We have an exact sequence
		\[
		k\longrightarrow U\longrightarrow UV
		\xrightarrow{q}Q\longrightarrow k.
		\]
		In particular, $(UV)^{\operatorname{co}q}=U$ and $q(V)=Q$.
		Moreover,
		\[
		V^{\operatorname{co}(q|_V)}
		=\{v\in V:(\operatorname{id}\otimes q)\Delta(v)=v\otimes1\}
		=U\cap V.
		\]
		Note that $U\cap V$ is a normal Hopf subalgebra of $V$. Then there is another exact sequence
		\[
		k\longrightarrow U\cap V\longrightarrow V
		\xrightarrow{q|_V}Q\longrightarrow k.
		\]
		By the Nichols--Zoeller theorem \cite{NZ}, we have
		$\dim(UV)=\dim U\dim Q$ and
		$\dim V=\dim(U\cap V)\dim Q$.
		
		Since $U$ and $V_0$ are cosemisimple Hopf algebras,
		by \cite[Corollary 5.1.12(a)]{RadfordBook}, $(UV)_0\subseteq U_0V_0=UV_0$.
		Furthermore, $S^2$ is the identity on $UV_0$. Observe that $UV_0$ is a Hopf subalgebra. Then $UV_0$ is cosemisimple and hence $UV_0\subseteq(UV)_0$. Consequently
		$(UV)_0=UV_0$.
	\end{proof}

	\section{Hopf algebras over a base scheme}\label{sec:families}
	
	We use the standard terminology and notation for schemes,
	geometric fibers, and finite locally free sheaves; see,
	e.g., \cite[Sections 10.3--10.5, 14.1, and 14.5]{Vakil}.
	Let $k$ be a field of characteristic zero and $X$ a
	$k$-scheme. A finite locally free Hopf algebra over $X$
	is a finite locally free $\mathcal O_X$-module equipped
	with $\mathcal O_X$-linear Hopf algebra structure maps.
	Let $\mathcal H$ be such a Hopf algebra of constant rank $N$.
	
	For a geometric point
	$\bar x:\Spec\Omega_{\bar x}\to X$ lying over $x\in X$, where
	$\Omega_{\bar x}$ is an algebraically closed field,
	the fiber $\mathcal H_{\bar x}:=\bar x^*\mathcal H$ is an $N$-dimensional Hopf algebra over $\Omega_{\bar x}$.
	Here pullbacks to $\Spec\Omega_{\bar x}$ are identified with
	their modules of global sections. If $\bar x$ factors through
	an affine open $U=\Spec R$, then
	$\mathcal H_{\bar x}\simeq\mathcal H(U)\otimes_R\Omega_{\bar x}$.
	
	Let
	\[
	\mathcal A=\mathcal H^\vee
	=\mathcal Hom_{\mathcal O_X}(\mathcal H,\mathcal O_X).
	\]
	Since $\mathcal H$ is finite locally free, $\mathcal A$
	is again a finite locally free Hopf algebra. Moreover,
	duals of finite locally free modules commute with pullback,
	as one checks locally on a basis. Hence, for every geometric
	point $\bar x$, we have
	\[
	\mathcal A_{\bar x}:=\bar x^*\mathcal A
	\simeq(\mathcal H_{\bar x})^\vee.
	\]
	
	We use the regular trace form of $\mathcal A$ to identify
	the Jacobson radical of each fiber. The following lemma
	describes the radical of this form. Let $E$ be a finite-dimensional algebra. For any $a\in E$, denote by 
	$L_a:E\to E$ the left multiplication operator given by $L_a(b)=ab$.
	
	\begin{lemma}\label{lem:trace}
		Let $E$ be a finite-dimensional algebra over an algebraically
		closed field of characteristic zero. The radical of the bilinear
		form
		\[
		T_E(a,b)=\Tr(L_{ab}:E\to E)
		\]
		is $J(E)$, where the radical is
		$\{a\in E:T_E(a,b)=0\text{ for all }b\in E\}$.
	\end{lemma}
	
	\begin{proof}
		If $a\in J(E)$, then $ab\in J(E)$ for every $b\in E$.
		Since $J(E)$ is nilpotent, $L_{ab}$ is nilpotent and
		$T_E(a,b)=0$.
		
		Let $0=E_0\subset E_1\subset\cdots\subset E_n=E$
		be a composition series of the regular left $E$-module.
		For each isomorphism class of simple $E$-modules, choose
		a representative $L$ with action map
		$\rho_L:E\to\End_k(L)$, and let $n_L$ be its
		multiplicity in this series.
		Since each $E_i$ is a left ideal, $L_{ab}$ induces
		on $E_i/E_{i-1}$ the operator $x+E_{i-1}\longmapsto abx+E_{i-1}$.
		Under an $E$-module isomorphism $E_i/E_{i-1}\simeq L$,
		this operator is $\rho_L(ab)=\rho_L(a)\rho_L(b)$.
		Additivity of trace along the composition series gives
		\[
		T_E(a,b)=\sum_L n_L
		\Tr\bigl(\rho_L(a)\rho_L(b)\bigr).
		\]
		Since $k$ is algebraically closed,
		\[
		E/J(E)\simeq\prod_L\End_k(L).
		\]
		Under this isomorphism, the bilinear form induced by $T_E$
		is
		\[
		((x_L)_L,(y_L)_L)
		\longmapsto
		\sum_L n_L\Tr(x_Ly_L).
		\]
		Each pairing $(x,y)\mapsto\Tr(xy)$ on $\End_k(L)$
		is nondegenerate, and $n_L\ne0$ in characteristic zero.
		Thus the induced form on $E/J(E)$ is nondegenerate,
		and the radical of $T_E$ is $J(E)$.
	\end{proof}
	
	\subsection{The trace map and multiplication}
	
	Using the canonical identification
	$\mathcal A^\vee\simeq\mathcal H$, define the
	$\mathcal O_X$-linear map $\tau:\mathcal A\to\mathcal H$
	by
	$\tau(a)(b)=\Tr(L_{ab}:\mathcal A\to\mathcal A)$,
	where $a,b$ are local sections of $\mathcal A$ and
	$L_{ab}$ denotes left multiplication by $ab$. Since the trace of an endomorphism of a finite locally
	free module commutes with pullback, $\tau_{\bar x}$
	is the regular trace map of $\mathcal A_{\bar x}$.
	From Lemma~\ref{lem:trace}, $\ker\tau_{\bar x}=J(\mathcal A_{\bar x})$.
	The trace form is symmetric. Then
	$\operatorname{im}\tau_{\bar x}=(\ker\tau_{\bar x})^\perp$ and so
	$J(\mathcal A_{\bar x})^\perp=(\mathcal H_{\bar x})_0$.
	As a result,
	\begin{equation}\label{eq:trace-coradical}
		\operatorname{im}\tau_{\bar x}
		=(\ker\tau_{\bar x})^\perp
		=J(\mathcal A_{\bar x})^\perp
		=(\mathcal H_{\bar x})_0.
	\end{equation}
	
	We first recall a dimension bound.
	
	\begin{lemma}\label{lem:family-generation-bound}
		Let $E$ be an $N$-dimensional unital algebra and 
		$C$ a linear subspace of $E$ containing $1$.
		For $r\ge1$, denote by $C^r$ the linear span of
		products of $r$ elements of $C$. Set
		$d=\max\{1,N-1\}$. Then $E$ is generated by $C$
		if and only if $C^d=E$.
	\end{lemma}
	
	\begin{proof}
		Since $1\in C$, we have $C^r\subseteq C^{r+1}$ and
		$\langle C\rangle_{\mathrm{alg}}=\bigcup_{r\ge1}C^r$.
		If $C^r=C^{r+1}$ for some $r$, multiplication by $C$
		shows that all subsequent powers equal $C^r$.
		Thus, if $C$ generates $E$, each inclusion is strict
		until the powers equal $E$.
		
		Suppose that $N>1$ and that $C$ generates $E$.
		Then $\dim C\ge2$. If $C^{N-1}\ne E$, all inclusions
		$C\subseteq C^2\subseteq\cdots\subseteq C^{N-1}$ are strict, so
		\[
		\dim C^{N-1}
		\ge \dim C+N-2
		\ge N,
		\]
		a contradiction. Thus $C^{N-1}=E$.
		If $N=1$, then $C=E$.
		The converse is immediate.
	\end{proof}
	
	For $r\ge1$, let
	$m_{\mathcal H}^{(r)}:\mathcal H^{\otimes r}\to\mathcal H$
	denote multiplication, with
	$m_{\mathcal H}^{(1)}=\operatorname{id}$.
	Put $d=\max\{1,N-1\}$ and consider the composite
	\begin{equation}\label{eq:family-generation}
		\mathcal G_d:
		\mathcal A^{\otimes d}
		\xrightarrow{\ \tau^{\otimes d}\ }
		\mathcal H^{\otimes d}
		\xrightarrow{\ m_{\mathcal H}^{(d)}\ }
		\mathcal H.
	\end{equation}
	
	Locally on $X$, choose bases of $\mathcal A^{\otimes d}$
	and $\mathcal H$. The matrix of $\mathcal G_d$ has
	$N^d$ columns and $N$ rows. By the Cauchy--Binet formula,
	the ideal generated by the $N\times N$ minors is unchanged
	under changes of bases. Denote the resulting ideal sheaf
	by $\mathcal I_N(\mathcal G_d)$.
	
	\begin{theorem}\label{thm:family}
		Let $k$ have characteristic zero, $X$ a $k$-scheme,
		and $\mathcal H$ a finite locally free Hopf algebra
		of constant rank $N$ over $X$. The locus 
		\[
		X^\circ=\{x\in X: B(\mathcal H_{\bar x})=\mathcal H_{\bar x}\text{ for a geometric point }\bar x\text{ over }x\}
		\]
		is open, and its complement is the zero locus of
		$\mathcal I_N(\mathcal G_d)$. Moreover, the ideal
		sheaf $\mathcal I_N(\mathcal G_d)$ and the open locus
		$X^\circ$ are compatible with arbitrary base change.
	\end{theorem}
	
	\begin{proof}
		For every geometric point $\bar x$, equation
		\eqref{eq:trace-coradical} shows that the image of
		$(\mathcal G_d)_{\bar x}$ is
		$(\mathcal H_{\bar x})_0^d$.
		Since $B(\mathcal H_{\bar x})$ is the algebra generated
		by $(\mathcal H_{\bar x})_0$, Lemma~\ref{lem:family-generation-bound}
		shows that $(\mathcal G_d)_{\bar x}$ is surjective if and
		only if $B(\mathcal H_{\bar x})=\mathcal H_{\bar x}$.
		Since $\mathcal H_{\bar x}$ has dimension $N$, this is
		equivalent to the nonvanishing of at least one
		$N\times N$ minor of $(\mathcal G_d)_{\bar x}$.
		Thus $X^\circ$ is open and its complement is the zero
		locus of $\mathcal I_N(\mathcal G_d)$.
		
		Let $f:Y\to X$ be a morphism. Since $\mathcal H$ is
		finite locally free, $f^*\mathcal A\simeq
		(f^*\mathcal H)^\vee$.
		The trace of an endomorphism of a finite locally free
		module, multiplication and tensor products commute
		with base change. Hence
		$f^*\mathcal G_d=\mathcal G_d(f^*\mathcal H)$ under
		the canonical identifications.
		The minors of the pulled-back map are the pullbacks
		of the corresponding minors of $\mathcal G_d$.
		Consequently, the ideal sheaf and the open locus
		have the asserted base change properties.
	\end{proof}

	\subsection{The matrix of \texorpdfstring{$\mathcal G_d$}{Gd}}
	On an open subset $W\subseteq X$ over which $\mathcal A$
	is free, choose an $\mathcal O_W$-basis
	$u_1,\ldots,u_N$ of $\mathcal A|_W$. Via the canonical
	isomorphism $\mathcal H|_W\simeq(\mathcal A|_W)^\vee$,
	let $f_1,\ldots,f_N$ be the dual basis. Then there are unique coefficients
	$m_{ij}^r,n_{ij}^r\in\mathcal O_X(W)$ such that
	\[
	u_i u_j=\sum_r m_{ij}^r u_r,\qquad
	f_i f_j=\sum_r n_{ij}^r f_r.
	\]
	Since the multiplication of $\mathcal A=\mathcal H^\vee$
	is dual to the coproduct of $\mathcal H$, one has
	\[
	\Delta(f_r)=\sum_{i,j}m_{ij}^r f_i\otimes f_j.
	\]
	Furthermore, the matrix of the regular trace form with respect to the  basis
	$u_1,\ldots,u_N$ has entries
	\[
	T_{ij}=\Tr(L_{u_i u_j})
	=\sum_{r,s}m_{is}^r m_{jr}^s.
	\]
	In particular, $\tau(u_i)=\sum_jT_{ij}f_j$.
	
	For $d\ge1$, let $n_{j_1,\ldots,j_d}^r$ be the
	coefficients determined by
	$f_{j_1}\cdots f_{j_d}
	=\sum_r n_{j_1,\ldots,j_d}^r f_r$.
	Associativity forces
	\[
	n_{j_1}^r=\delta_{j_1r},\qquad
	n_{j_1,\ldots,j_d}^r
	=\sum_a n_{j_1,\ldots,j_{d-1}}^a n_{a j_d}^r
	\quad(d\ge2).
	\]
	With respect to the bases
	$u_{i_1}\otimes\cdots\otimes u_{i_d}$ and $f_r$,
	the matrix of $\mathcal G_d$ has entries
	\begin{equation}\label{eq:general-parameter-matrix}
		(\mathcal G_d)_{r,(i_1,\ldots,i_d)}
		=\sum_{j_1,\ldots,j_d}
		T_{i_1j_1}\cdots T_{i_dj_d}
		n_{j_1,\ldots,j_d}^r.
	\end{equation}
	In particular, every $N\times N$ minor is a polynomial
	in the structure constants $m_{ij}^r$ and $n_{ij}^r$.
	
	When $f_1,\ldots,f_N$ is a fixed PBW basis of
	$\mathcal H$, the coefficients $n_{ij}^r$ are obtained
	by multiplying PBW monomials and using the defining
	relations, while the coefficients $m_{ij}^r$ are
	determined by the coproduct. If these structure constants
	are polynomial in the parameters of a family of liftings,
	every maximal minor in
	\eqref{eq:general-parameter-matrix} is polynomial in
	those parameters. More generally, if the structure
	constants are regular on an open subset, the minors
	are regular there.
	
	We also record the change of basis. Suppose that
	$f'_1,\ldots,f'_N$ is another basis of $\mathcal H$. Then there exists $U\in\GL_N(\mathcal O_X(W))$ such that
	$(f'_1,\ldots,f'_N)=(f_1,\ldots,f_N)U$. The corresponding dual
	basis satisfies
	$(u'_1,\ldots,u'_N)=(u_1,\ldots,u_N)U^{-\mathsf t}$.
	Consequently, the matrices of $\mathcal G_d$ satisfy
	\begin{equation}\label{eq:general-chart-change}
		\mathcal G'_d
		=U^{-1}\mathcal G_d(U^{-\mathsf t})^{\otimes d}.
	\end{equation}
	By the Cauchy--Binet formula, every
	$N\times N$ minor of $\mathcal G'_d$ belongs to the
	ideal generated by the $N\times N$ minors of
	$\mathcal G_d$. The reverse inclusion follows by the inverse change of basis. Thus the ideals generated by the two sets of
	minors are the same.
	
	\begin{corollary}\label{cor:generic}
		Suppose that $Y$ is a reduced irreducible component of $X$, and that $\eta$ is its generic point. Set
		$Y^\circ=Y\cap X^\circ$, and put
		$\mathfrak b=\mathcal I_N(\mathcal G_d)\mathcal O_Y$.
		The following conditions are equivalent:
		\begin{enumerate}
			\item $Y^\circ\ne\varnothing$;
			\item $\mathcal G_d$ has rank $N$ at $\eta$;
			\item $\mathfrak b\ne0$.
		\end{enumerate}
		If these conditions hold, $Y^\circ$ is dense open in $Y$.
	\end{corollary}
	
	\begin{proof}
		By Theorem~\ref{thm:family},  $Y^\circ$ is open and its complement is the zero locus of
		$\mathfrak b$. Since $Y$ is irreducible, $Y^\circ$ is nonempty if and only if it contains the generic point $\eta$.
		
		On a nonempty affine open subset of $Y$ where
		$\mathcal A^{\otimes d}$ and $\mathcal H$ are free,
		the coordinate ring embeds into the function field
		$\kappa(\eta)$, since $Y$ is reduced and irreducible.
		Hence $\mathfrak b\ne0$ if and only if at least one
		$N\times N$ minor is nonzero at $\eta$. Consequently, the conditions are equivalent, and then the nonempty open subset
		$Y^\circ$ is dense in $Y$.
	\end{proof}
	
	\begin{remark}
		The open subset $X^\circ$ may be empty. For example,
		let $A=\mathcal B(V)\#k\Gamma$ with $V\ne0$ and $\dim A<\infty$.
		The kernel of the canonical Hopf algebra projection $A\to k\Gamma$
		is the nonzero nilpotent Hopf ideal
		$I=\mathcal B(V)^+\#k\Gamma$, so $I\subseteq J(A)$.
		Since $I$ is a Hopf ideal, $I^\perp$ is a proper
		Hopf subalgebra of $A^*$. By
		Lemma~\ref{lem:dual-radical-annihilator},
		$(A^*)_0=J(A)^\perp\subseteq I^\perp$.
		Hence $B(A^*)\subseteq I^\perp\ne A^*$.
		Thus the constant family with fiber $A^*$
		has $X^\circ=\varnothing$.
	\end{remark}

	\section{Hopf ideals in liftings of Nichols algebras}\label{sec:nichols}
	
	In this section, $k$ is algebraically closed of characteristic zero. Let $A$ be a finite-dimensional pointed Hopf algebra over $k$ with abelian coradical $A_0=k\Gamma$. Then by \cite[Theorem 4.15]{AngionoRelations}, $\gr A\simeq\B(V)\#k\Gamma$, where $V\in\YD$  is the infinitesimal braiding of $A$ and is of diagonal type. Namely, $A$ is a lifting of a finite-dimensional Nichols algebra of diagonal type.
	
	In Section~\ref{sec:families}, the equality
	$B(A^*)=A^*$ is characterized by the maximal minors
	of $\mathcal G_d$, whose entries are determined by the
	multiplication and coproduct of $A$.
	
	Let $E_{g,\chi}=\iota(V_g^\chi)$, where
	$\iota:V\to A_1$ is introduced in Section~\ref{sec:prelim}.
	Every character occurring in $V$ is nontrivial.
	We prove that every Hopf ideal $K$ of $A$ with
	$K\cap k\Gamma=0$ is generated by its intersections
	with the spaces $E_{g,\chi}$.
	If $R=\bigoplus_{g,\chi}(J(A)\cap E_{g,\chi})$,
	then
	$\Imax(A)=(R)$ and $B(A^*)=(A/(R))^*$.
	Now we determine $J(A)\cap E_{g,\chi}$ by using the defining relations of $A$.

	\subsection{Hopf ideals and skew-primitive generators}
	
	We begin with the Hopf ideals of $A$ that meet
	$k\Gamma$ trivially. Let $K$ be a Hopf ideal of $A$ with
	$K\cap k\Gamma=0$, and let $\pi:A\to A/K$ be the quotient
	map. Equip $A/K$ with the quotient filtration
	$\overline A_n=\pi(A_n)=(A_n+K)/K$.
	
	The map $\pi$ can induce a surjective graded Hopf map
	\[
	\operatorname{gr}\pi:\operatorname{gr}A
	\longrightarrow\operatorname{gr}(A/K),
	\qquad
	a+A_{n-1}\longmapsto
	\pi(a)+\overline A_{n-1}
	\quad(a\in A_n).
	\]
	Let
	$\theta:\mathcal B(V)\#k\Gamma\to\operatorname{gr}A$
	be the graded Hopf algebra isomorphism induced by
	\[
	1\#g\longmapsto g+A_{-1},\qquad
	v\#1\longmapsto\iota(v)+A_0
	\]
	for $g\in\Gamma$ and $v\in V$. Set
	$\phi_K=\operatorname{gr}\pi\circ\theta$. Then
	\[
	\phi_K(1\#g)=\pi(g),\qquad
	\phi_K(v\#1)=\pi(\iota(v))+\overline A_0.
	\]
	In particular,
	\[
	\phi_K(v_1\cdots v_r\#g)
	=
	\bigl(\pi(\iota(v_1))+\overline A_0\bigr)
	\cdots
	\bigl(\pi(\iota(v_r))+\overline A_0\bigr)\pi(g).
	\]
	Since $K\cap k\Gamma=0$, $\phi_K$ is an isomorphism
	in degree zero.

	\begin{lemma}\label{lem:degree-one-kernel}
		There exists a subobject $R'$ of $V$ in $\YD$ such that $(\ker\phi_K)_1=k\Gamma R'$.
	\end{lemma}
	
	\begin{proof}
		Let $\mathcal C=\mathcal B(V)\#k\Gamma$. Then
		$\mathcal C(1)=V\#k\Gamma$ is a $k\Gamma$-bicomodule
		with coactions
		\[
		\lambda(v\#h)=gh\otimes(v\#h),\qquad
		\rho(v\#h)=(v\#h)\otimes h
		\]
		for $v\in V_g$ and $h\in\Gamma$.
		Since $\phi_K$ is a graded Hopf algebra map and is
		injective in degree zero, $W=(\ker\phi_K)_1$ is a
		sub-bicomodule of $\mathcal C(1)$.
		Consequently,
		\[
		W=\bigoplus_{g,h\in\Gamma}
		\bigl(W\cap(V_g\#h)\bigr).
		\]
		
		Set $R'=\{v\in V:v\#1\in W\}$.
		Since $\ker\phi_K$ is a two-sided ideal of $\mathcal C$,
		right multiplication by $1\#h$ and $1\#h^{-1}$ preserves
		$W$. Hence
		\[
		W\cap(V_g\#h)
		=(R'\cap V_g)\#h.
		\]
		Thus $W=R'\#k\Gamma$ and
		$R'=\bigoplus_g(R'\cap V_g)$.
		
		For $h\in\Gamma$ and $v\in R'$, we also have
		\[
		(h\cdot v)\#1
		=(1\#h)(v\#1)(1\#h^{-1})\in W.
		\]
		Hence $R'$ is stable under the $\Gamma$-action and so
		$R'$ is a subobject of $V$ in $\YD$.
		Therefore $W=R'\#k\Gamma=k\Gamma R'$.
	\end{proof}

	\begin{lemma}\label{lem:degree-one-identification}
		Let $R'$ be the subobject of $V$ in
		Lemma~\ref{lem:degree-one-kernel}, and $R_K=\bigoplus_{g,\chi}(K\cap E_{g,\chi})$. 
		Then $\iota(R')=R_K$ and $K\cap A_1=k\Gamma R_K$.
	\end{lemma}
	
	\begin{proof}
		Since $R'$ is a subobject of $V$ in $\YD$, we have
		\[
		R'=\bigoplus_{g,\chi}(R'\cap V_g^\chi).
		\]
		Let $v\in R'\cap V_g^\chi$. By
		Lemma~\ref{lem:degree-one-kernel},
		$v\#1\in(\ker\phi_K)_1$, and hence $\pi(\iota(v))\in\pi(A_0)$.
		For $h\in\Gamma$, we have
		\[
		\pi(h)\pi(\iota(v))\pi(h^{-1})=\chi(h)\pi(\iota(v)).
		\]
		Since $\Gamma$ is abelian, conjugation by $\Gamma$
		acts trivially on $\pi(A_0)$. Thus
		$\pi(\iota(v))=0$, as $\chi\ne\varepsilon$.
		Consequently, $\iota(v)\in K\cap E_{g,\chi}$, and
		$\iota(R')\subseteq R_K$.
		
		Conversely, let $v\in V_g^\chi$ with
		$\iota(v)\in K$. Then $\phi_K(v\#1)=0$, so
		$v\#1\in(\ker\phi_K)_1=k\Gamma R'$.
		As a result, $v\in R'$ and so $\iota(R')=R_K$.
		
		Now let $x\in K\cap A_1$. Its class in $A_1/A_0$,
		identified with $\mathcal C(1)$, belongs to
		$(\ker\phi_K)_1=k\Gamma R'$. There exists
		$y\in k\Gamma R_K$ such that $x-y\in A_0$.
		Since $K$ is a two-sided ideal and $R_K\subseteq K$,
		we have $y\in K$. Therefore
		$x-y\in K\cap A_0=0$, and consequently
		$K\cap A_1=k\Gamma R_K$.
	\end{proof}

	\begin{lemma}\label{lem:graded-quotient-first-layer}
		Let $R'$ be the subobject of $V$ in
		Lemma~\ref{lem:degree-one-kernel} and $U\in\YD$ such
		that $V=R'\oplus U$.
		Then $\mathcal B(U)\#k\Gamma\cong\gr(A/K)$ as
		graded Hopf algebras.
	\end{lemma}
	
	\begin{proof}
		Let $j:U\to V$ be the inclusion and
		$p:V\to U$ the projection with kernel $R'$.
		By \cite[Corollary 2.3]{ASPointed}, the maps $j$ and
		$p$ induce morphisms of Nichols algebras
		\[
		\B(U)\xrightarrow{\B(j)}\B(V)\xrightarrow{\B(p)}\B(U),
		\]
		with $\B(p)\B(j)=\id_{\B(U)}$. Let $q:\B(V)\to\B(V)/(R')$ be the quotient map.
		Since $\B(p)$ vanishes on $R'$, the map $\B(p)$
		factors through $q$. Since $\B(V)/(R')$ is generated by the image of $U$, $q\B(j)$ is surjective.
		Moreover, $\overline{\B(p)}q\B(j)=\id_{\B(U)}$.
		Consequently, $q\B(j)$ is an isomorphism and
		$\B(V)/(R')\cong\B(U)$.
		
		By Lemma~\ref{lem:degree-one-kernel},
		$R'\#1\subseteq\ker\phi_K$. Since $\ker\phi_K$
		is a Hopf ideal, $\ker\phi_K$ contains the Hopf ideal
		generated by $R'\#1$, namely $(R')\#k\Gamma$.
		Thus $\phi_K$ induces a surjective graded Hopf map
		$\overline{\phi}_K:
		(\B(V)/(R'))\#k\Gamma\to\gr(A/K)$.
		The isomorphism $q\B(j)$ induces a graded Hopf
		algebra isomorphism
		$(q\B(j))\#\id_{k\Gamma}$ between the corresponding
		bosonizations. Therefore, the composition
		\[
		\psi_K=
		\overline{\phi}_K\circ
		\bigl((q\B(j))\#\id_{k\Gamma}\bigr):
		\B(U)\#k\Gamma\longrightarrow\gr(A/K)
		\]
		is an epimorphism of graded Hopf algebras.
		
		The restriction of $\psi_K$ to the degree zero component
		is injective, since $\phi_K$ is injective in degree zero.
		The restriction of $\psi_K$ to the degree-one component
		$U\#k\Gamma$ agrees with the restriction of $\phi_K$
		to $U\#k\Gamma$. By Lemma~\ref{lem:degree-one-kernel},
		$(\ker\phi_K)_1=R'\#k\Gamma$. Since $V=R'\oplus U$,
		we have $(U\#k\Gamma)\cap(R'\#k\Gamma)=0$.
		Hence $\psi_K$ is injective in degree one.
		
		The bosonization $C=\B(U)\#k\Gamma$ is
		coradically graded, with $C_n=\bigoplus_{j=0}^n C(j)$.
		In particular, $C_1=C(0)\oplus C(1)$.
		Then $\psi_K\mid_{C_1}$ is injective.
		By \cite[Theorem 5.3.1]{Montgomery},  $\psi_K$ is an isomorphism.
	\end{proof}
	
	\begin{remark}
		The isomorphism class in Lemma~\ref{lem:graded-quotient-first-layer}
		does not depend on the choice of $U$. Let $U'$ be another Yetter--Drinfeld complement of
		$R'$ in $V$, and let $p':V\to U'$ be the projection
		with kernel $R'$. The restriction $p'|_U:U\to U'$
		is an isomorphism in $\YD$. By the functoriality
		of Nichols algebras, $\B(p'|_U)$ induces an
		isomorphism $\B(U)\#k\Gamma\cong\B(U')\#k\Gamma$
		of graded Hopf algebras. Moreover, the isomorphisms
		with $\gr(A/K)$ constructed in the proof are compatible
		with $\B(p'|_U)\#\id_{k\Gamma}$. 
	\end{remark}
	
	\begin{theorem}\label{thm:first-layer}
		Let $K$ be a Hopf ideal of $A$ with $K\cap k\Gamma=0$ and
		put $R_K=\bigoplus_{g,\chi}(K\cap E_{g,\chi})$.
		Then $K=(R_K)$.
	\end{theorem}
	
	\begin{proof}
		Let $I=(R_K)$. Every $v\in R_K\cap E_{g,\chi}$
		is $(g,1)$-skew primitive, with $\varepsilon(v)=0$
		and $S(v)=-g^{-1}v$. Hence $I$ is a Hopf ideal.
		Since $I\subseteq K$ and $R_K\subseteq I$, we have
		$I\cap E_{g,\chi}=K\cap E_{g,\chi}$ for all $g,\chi$.
		In particular, $R_I=R_K$.
		
		By Lemma~\ref{lem:degree-one-identification}, the
		Yetter--Drinfeld submodules associated with $I$ and
		$K$ both equal $\iota^{-1}(R_K)$. Choose a
		Yetter--Drinfeld complement $U$ of $\iota^{-1}(R_K)$
		in $V$. Lemma~\ref{lem:graded-quotient-first-layer}
		applies to both $I$ and $K$, and yields isomorphisms
		\[
		\B(U)\#k\Gamma
		\xrightarrow{\ \sim\ }\gr(A/I),
		\qquad
		\B(U)\#k\Gamma
		\xrightarrow{\ \sim\ }\gr(A/K).
		\]
		
		Let $\pi_I:A\to A/I$ and $\pi_K:A\to A/K$
		be the quotient maps, and let
		$\overline\pi:A/I\to A/K$ be the map induced by
		$I\subseteq K$. From $\pi_K=\overline\pi\pi_I$, one has
		$\phi_K=\gr(\overline\pi)\phi_I$.
		Since the two isomorphisms above are obtained using
		the same complement $U$ and the same quotient
		$\B(V)/( \iota^{-1}(R_K))$, the induced maps satisfy
		\[
		\gr(\overline\pi)\circ\psi_I=\psi_K.
		\]
		Consequently, $\gr(\overline\pi)$ is an isomorphism.
		
		As $A$ is finite-dimensional,
		$\dim(A/I)=\dim\gr(A/I)=\dim\gr(A/K)=\dim(A/K)$.
		Thus $\dim I=\dim K$ and consequently $I=K$.
	\end{proof}
	
	\begin{corollary}\label{cor:general-largest}
		Let $R=\bigoplus_{g,\chi}\bigl(J(A)\cap E_{g,\chi}\bigr)$. 
		Then $\Imax(A)=(R)$ and $B(A^*)=(A/(R))^*$.
	\end{corollary}
	\begin{proof}
		Clearly, $(R)$ is a Hopf ideal contained in  $J(A)$ and so $(R)\subseteq\Imax(A)$.
		Since $k\Gamma$ is semisimple, $\Imax(A)\cap k\Gamma\subset J(A)\cap k\Gamma=0$.
		
		Let $R_{\Imax(A)}=\bigoplus_{g,\chi}(\Imax(A)\cap E_{g,\chi})$.	By Theorem~\ref{thm:first-layer},
		$\Imax(A)=(R_{\Imax(A)})$. Since
		$\Imax(A)\subseteq J(A)$, we have
		$R_{\Imax(A)}\subseteq R$, and consequently
		$\Imax(A)\subseteq(R)$.
		Therefore $\Imax(A)=(R)$.
		
		By Lemma~\ref{lem:dual-radical-annihilator}, $\Imax(A)=B(A^*)^\perp$ and so
		$B(A^*)=\Imax(A)^\perp=(R)^\perp=(A/(R))^*$.
	\end{proof}

	\subsection{The action on the defining relations}
	
	By Theorem~\ref{thm:first-layer} and
	Corollary~\ref{cor:general-largest}, it remains to
	determine $J(A)\cap E_{g,\chi}$ for each pair $(g,\chi)$.
	We use the defining relations of $A$ to characterize
	these intersections.
	
	Let $\mathcal G=\Aut_{\YD}(V)$, regarded as the algebraic
	subgroup of $\GL(V)$ consisting of automorphisms of $V$
	in $\YD$. Since $\Gamma$ is finite abelian and $k$ is
	algebraically closed, we have $V=\bigoplus_{g,\chi}V_g^\chi$.
	Every automorphism in $\mathcal G$ preserves the
	$\Gamma$-action and coaction, and hence each
	$V_g^\chi$. Conversely, any collection of invertible
	linear maps on the components $V_g^\chi$ defines
	an automorphism in $\YD$. Thus, as algebraic groups,
	\[
	\mathcal G=\prod_{g,\chi}\GL(V_g^\chi).
	\]
	See \cite{ASPointed} for more details.
	
	By \cite[Corollary 2.3]{ASPointed}, every
	$F\in\mathcal G$ induces a graded Hopf algebra
	automorphism of $\B(V)$. Hence the induced action
	of $\mathcal G$ on $T(V)$ preserves $\mathcal J(V)$,
	and each homogeneous component $\mathcal J^m(V)$
	is $\mathcal G$-stable.
	
	Let $W$ be a finite-dimensional homogeneous $\mathcal G$-stable
	subspace of $\mathcal J(V)$ such that $\mathcal J(V)=(W)$.
	For $m\ge2$, put $W_m=W\cap V^{\otimes m}$. Then
	\[
	W=\bigoplus_{m\ge2}W_m,
	\]
	with only finitely many nonzero summands.
	
	Suppose that the fixed lifting $A$ has a presentation
	\begin{equation}\label{eq:relation-presentation}
		A=(T(V)\#k\Gamma)/(r-\kappa(r):r\in W),
		\qquad \kappa\in\Hom_k(W,k\Gamma).
	\end{equation}
	Let $\pi:T(V)\#k\Gamma\to A$ be the quotient map.
	The restriction of $\pi$ to $V$ is $\iota$.  
	
	Fix $E=V_g^\chi$. For $v\in E$ and $\varphi\in E^*$,
	extend $\varphi$ to $V$ by setting $\varphi$ equal
	to zero on every $V_{g'}^{\chi'}$ with
	$(g',\chi')\ne(g,\chi)$.
	Let $\delta_{v,\varphi}:T(V)\to T(V)$ be the
	derivation determined by
	$\delta_{v,\varphi}(w)=\varphi(w)v$ for $w\in V$.
	Explicitly,
	\[
	\delta_{v,\varphi}(w_1\cdots w_m)
	=\sum_{j=1}^m\varphi(w_j)
	w_1\cdots w_{j-1}v w_{j+1}\cdots w_m.
	\]
	
	Under the identification
	\[
	\Lie(\mathcal G)
	=\bigoplus_{g,\chi}\End_k(V_g^\chi),
	\]
	the endomorphism $v\otimes\varphi$ belongs to
	$\Lie(\mathcal G)$. The differential of the
	$\mathcal G$-action on $V^{\otimes m}$ is
	$\delta_{v,\varphi}$. Since $W$ is
	$\mathcal G$-stable, we have
	$\delta_{v,\varphi}(W)\subseteq W$.
	
	Define
	\begin{equation}\label{eq:intrinsic-contraction}
		\Psi_E:E\longrightarrow\Hom_k(E^*\otimes W,k\Gamma),
		\qquad
		\Psi_E(v)(\varphi\otimes r)
		=\kappa(\delta_{v,\varphi}r).
	\end{equation}

	\begin{theorem}\label{thm:contraction}
		With $W$ as above and $A$ presented as in
		\eqref{eq:relation-presentation},
		we have
		\[
		\ker\Psi_E=E\cap\iota^{-1}(J(A)).
		\]
	\end{theorem}

	\begin{proof}
		Suppose that $v\in E$ and $\iota(v)\in J(A)$. Then every term of $\delta_{v,\varphi}r$ contains $v$ and so
		$\pi(\delta_{v,\varphi}r)\in J(A)$.
		On the other hand, $\delta_{v,\varphi}r\in W$, which implies that $\pi(\delta_{v,\varphi}r)=\kappa(\delta_{v,\varphi}r)\in k\Gamma$.
		Since $k\Gamma$ is semisimple, $J(A)\cap k\Gamma=0$.
		Hence
		$\kappa(\delta_{v,\varphi}r)=0$ for every
		$\varphi\in E^*$ and $r\in W$. Consequently,
		$v\in\ker\Psi_E$.
		
		Conversely, let $0\ne v\in\ker\Psi_E$.
		Choose $\varphi\in E^*$ with $\varphi(v)=1$,
		and extend $\varphi$ by zero on the other
		isotypic components.
		Then $U=\ker\varphi$ is a subobject of $V$ in $\YD$ such that
		$V=kv\oplus U$. Furthermore, there is an algebraic group homomorphism
		\[
		\lambda:\mathbb G_m\longrightarrow\mathcal G,
		\qquad
		\lambda(t)=t\id_{kv}\oplus\id_U.
		\]
		
		The homomorphism $\lambda$ induces an action of
		$\mathbb G_m$ on $T(V)$ by algebra automorphisms.
		For a tensor word $w_1\cdots w_m$, with each
		$w_j\in kv$ or $w_j\in U$, the action is
		\[
		\lambda(t)\cdot(w_1\cdots w_m)
		=t^q w_1\cdots w_m,
		\]
		where $q$ is the number of factors belonging to $kv$.
		Since $W$ is $\mathcal G$-stable, the action of
		$\mathbb G_m$ on $T(V)$ preserves $W$. The
		restriction to $W$ is a finite-dimensional
		rational representation of $\mathbb G_m$, and hence
		\[
		W=\bigoplus_{q\ge0}W^{(q)},
		\]
		where $W^{(q)}$ consists of the elements $r\in W$
		satisfying $\lambda(t)\cdot r=t^q r$ for all
		$t\in\mathbb G_m$.
		
		The differential at $t=1$ of the
		$\mathbb G_m$-action on $T(V)$ is
		$\delta_{v,\varphi}$. Thus
		$\delta_{v,\varphi}$ acts on $W^{(q)}$
		as multiplication by $q$. For $r\in W^{(q)}$,
		$\Psi_E(v)=0$ and hence
		\[
		q\kappa(r)
		=\kappa(\delta_{v,\varphi}r)=0.
		\]
		Since $k$ has characteristic zero,
		$\kappa(r)=0$ whenever $q>0$.
		
		Assign degree one to $v$ and degree zero to
		$U$ and $\Gamma$. Every element of $W^{(q)}$
		is homogeneous of degree $q$. For $r\in W^{(0)}$,
		$r-\kappa(r)$ is homogeneous of degree zero,
		whereas for $r\in W^{(q)}$ with $q>0$,
		$r-\kappa(r)=r$ is homogeneous of degree $q$.
		Since $W=\bigoplus_{q\ge0}W^{(q)}$ and $\kappa$
		is linear, the defining ideal in
		\eqref{eq:relation-presentation} is generated
		by homogeneous elements. Thus
		\eqref{eq:relation-presentation} defines a
		nonnegative grading of $A$.
		
		The coproduct and antipode preserve the grading.
		Indeed, for $x\in U\cap V_h^\psi$, and also
		for $x=v$ with $h=g$ and $\psi=\chi$,  
		\[
		\Delta(\iota(x))
		=\iota(x)\otimes1+h\otimes\iota(x),
		\qquad
		S(\iota(x))=-h^{-1}\iota(x).
		\]
		The elements of $\iota(U)$ and $k\Gamma$ have
		degree zero, whereas $\iota(v)$ has degree one.
		Since $U$ is a subobject of $V$ in $\YD$,
		every element of $U$ is a sum of elements in
		$U\cap V_h^\psi$. Therefore, the coproduct and antipode preserve the
		grading on the algebra generators. Thus $A$ is
		a graded Hopf algebra.
		
		Put $A_+=\bigoplus_{n>0}A(n)$, the ideal of
		positive-degree elements. There is $n>0$ such that
		$A(m)=0$ for all $m>n$ because $A$ is finite-dimensional.  Thus
		$A_+\subseteq J(A)$, and in particular,
		$\iota(v)\in J(A)$.
	\end{proof}

	\begin{corollary}\label{cor:contraction}
		Under the hypotheses of Theorem~\ref{thm:contraction},
		we have
		\[
		\Imax(A)=
		\left(\iota\left(\bigoplus_E\ker\Psi_E\right)\right),
		\]
		where the sum ranges over the nonzero isotypic
		components $E=V_g^\chi$. $A^*$ is
		generated by the coradical if and only if
		$\Psi_E$ is injective for every nonzero isotypic
		component $E=V_g^\chi$.
	\end{corollary}
	
	\begin{proof}
		By Theorem~\ref{thm:contraction} and
		Corollary~\ref{cor:general-largest},
		\[
		\Imax(A)=\left(
		\bigoplus_{g,\chi}
		J(A)\cap\iota(V_g^\chi)
		\right)\\
		=\left(
		\iota\left(
		\bigoplus_{g,\chi}\ker\Psi_{V_g^\chi}
		\right)
		\right).
		\]
		From Lemma~\ref{lem:dual-radical-annihilator},
		$A^*=B(A^*)$ is equivalent to
		$\Imax(A)=0$. Since $\iota$ is injective,  $\Imax(A)=0$ precisely when $\ker\Psi_E=0$ for every nonzero isotypic
		component $E=V_g^\chi$.
	\end{proof}
	
	For fixed bases of $E$, $W$ and $k\Gamma$, the matrix entries
	of $\Psi_E$ are linear in the coefficients of $\kappa$.
	If these coefficients are polynomial in the deformation parameters,
	then $\Psi_E$ is injective if and only if
	$\operatorname{rank}\Psi_E=\dim E$.

	\begin{remark}
		The description of $\Imax(A)$ in
		Corollary~\ref{cor:general-largest} holds for every
		finite-dimensional pointed Hopf algebra considered
		in this section. However, Theorem~\ref{thm:contraction} assumes that the defining relations of
		$\B(V)$ admit deformations of the form
		$r-\kappa(r)$, with $\kappa(r)\in k\Gamma$.
		In the recursive construction of liftings, deformations
		of the defining relations may involve terms
		outside $k\Gamma$; see \cite{AngionoSurvey} and reference therein.
	\end{remark}

	\section{Parameter criteria for liftings}\label{sec:parameters}
	
	Let $k$ be algebraically closed of characteristic zero and $A$ a finite-dimensional pointed Hopf algebra with abelian coradical. $A^*$ is generated by its coradical if and only if $J(A)\cap E_{g,\chi}=0$ for every Yetter--Drinfeld pair $(g,\chi)$; see Corollary~\ref{cor:general-largest}. We describe these intersections in terms of deformation parameters for liftings of quantum linear spaces and standard liftings of finite Cartan type.
	
	\subsection{Liftings of quantum linear spaces}\label{sec:qls}

	Let $\Gamma$ be a finite abelian group. Fix Yetter--Drinfeld pairs
	$(g_i,\chi_i)$, $1\le i\le n$, where
	$g_i\in\Gamma$ and $\chi_i\in\widehat\Gamma$. Put
	\[
	q_{ij}=\chi_j(g_i),\qquad
	N_i=\operatorname{ord}q_{ii}>1.
	\]
	
	\begin{definition}\cite[Lemma 3.4]{ASLifting}
		Let $V=\bigoplus_{i=1}^n kx_i\in\YD$, with
		$x_i\in V_{g_i}^{\chi_i}$ such that $q_{ij}q_{ji}=1$ for $i\ne j$. The Nichols algebra $\B(V)$ is called a quantum linear
		space, where the defining relations are given by 
		\[
		x_i^{N_i}=0,\qquad
		x_ix_j-q_{ij}x_jx_i=0\quad(i<j).
		\]
	\end{definition}

	Put
	\[
	\mathcal W_{\mathrm{qls}}
	=
	\operatorname{span}_k
	\left(
	\{x_i^{N_i}\mid 1\le i\le n\}
	\cup
	\{x_ix_j-q_{ij}x_jx_i\mid 1\le i<j\le n\}
	\right).
	\]
	
	\begin{definition}
		The collection $D=\bigl((a_i)_{1\le i\le n},(c_{ij})_{1\le i<j\le n}\bigr)$
		of elements in $k$ is called a compatible datum if
		\[
		\begin{aligned}
			a_i&=0
			&&\text{if }g_i^{N_i}=1
			\text{ or }\chi_i^{N_i}\ne\varepsilon,\\
			c_{ij}&=0
			&&\text{if }g_ig_j=1
			\text{ or }\chi_i\chi_j\ne\varepsilon.
		\end{aligned}
		\]
		For a compatible datum $D$, let $\kappa_D:\mathcal W_{\mathrm{qls}}\longrightarrow k\Gamma$
		be the linear map determined by
		\[
		\begin{aligned}
			\kappa_D(x_i^{N_i})
			&=a_i(1-g_i^{N_i}),\\
			\kappa_D(x_ix_j-q_{ij}x_jx_i)
			&=c_{ij}(1-g_ig_j),\qquad i<j.
		\end{aligned}
		\]
	\end{definition}
	
	For $i<j$, set $c_{ji}=-q_{ji}c_{ij}$.
	The power parameters in \cite{ASLifting} are normalized to zero or one;
	allowing arbitrary $a_i\in k$ amounts to rescaling the generators.
	By \cite[Lemma 5.1, Proposition 5.2 and Theorem 5.5]{ASLifting},
	every lifting of $\mathcal B(V)$ is isomorphic to the Hopf algebra
	\[
	A_D=(T(V)\#k\Gamma)/(r-\kappa_D(r):r\in\mathcal W_{\mathrm{qls}})
	\]
	for some compatible datum $D$. For every such $D$, the elements
	\[
	hx_1^{e_1}\cdots x_n^{e_n},\qquad
	h\in\Gamma,\quad0\le e_i<N_i,
	\]
	form a basis of $A_D$.
	
	\subsubsection{The largest Hopf ideal contained in the radical}
	\label{subsec:qls-radical}
	
	Throughout this subsection, fix a compatible datum $D$ and set $A=A_D$ for short.
	
	For a Yetter--Drinfeld pair $(g,\chi)$, set $E=E_{g,\chi}$. The
	corresponding set of indices is $I_E=\{i:(g_i,\chi_i)=(g,\chi)\}$ and $E=\operatorname{span}\{x_i:i\in I_E\}$.
	Write $N=\operatorname{ord}\chi(g)$. If $N>2$, then
	$E$ is one-dimensional; indeed, if $i,j\in I_E$ with
	$i\ne j$, then $q_{ij}=q_{ji}=\chi(g)$, whereas
	$q_{ij}q_{ji}=1$, which implies that $N\le2$.
	
	For $N>2$, write $E=kx_i$ and define
	\[
	R_E=
	\begin{cases}
		kx_i,&a_i=0\text{ and }c_{ij}=0
		\text{ for every }j\ne i,\\
		0,&\text{otherwise}.
	\end{cases}
	\]
	
	Suppose now that $N=2$. Define
	\[
	b_{ii}=2a_i,\qquad b_{ij}=c_{ij}\quad(i\ne j).
	\]
	Consider the linear map $	\Theta_E:E\longrightarrow\bigoplus_{j=1}^n k\Gamma$ given by 
	\begin{equation}\label{eq:pairing-map}
		\sum_{i\in I_E}t_ix_i
		\longmapsto
		\left(\left(\sum_{i\in I_E}t_ib_{ij}\right)(1-gg_j)\right)_{j=1}^n.
	\end{equation}
	Define $R_E=\ker\Theta_E$. Equivalently, $R_E$ is the
	kernel of the matrix over $k$ $(b_{ij})_{1\le j\le n,\,
		i\in I_E}$ after deleting the rows indexed by $j$ with
	$gg_j=1$. This follows from the linear independence of
	the group elements in $k\Gamma$ and the parameter
	conditions imposed on $D$.
	
	The subspace $R=\bigoplus_E R_E$ 
	is defined by these intersections, where $E$ ranges
	over the nonzero isotypic components $E_{g,\chi}$.
	
	\begin{lemma}\label{lem:qls-radical-intersections}
		For every nonzero isotypic component $E=E_{g,\chi}$, $J(A)\cap E=R_E$.
	\end{lemma}
	
	\begin{proof}
		Recall that $J(A)\cap k\Gamma=0$. Let $0\ne v\in R_E$. If $N>2$, then $E=kx_i$. By definition of $R_E$,
		$a_i=0$ and $c_{ij}=0$ for every $j\ne i$.
		Consequently,
		\[
		x_i^{N_i}=0,\qquad
		x_ix_j=q_{ij}x_jx_i\quad(j\ne i).
		\]
		The relations with $\Gamma$ also show that $x_i$
		commutes with every generator of $A$ up to multiplication
		by an element in $k^\times$. The ideal generated by $x_i$ is
		nilpotent, and so $x_i\in J(A)$.
		
		Now suppose that $N=2$ and write
		$v=\sum_{i\in I_E}t_ix_i$. For every
		$1\le j\le n$,
		\begin{equation}\label{eq:normal-radical}
			vx_j-\chi_j(g)x_jv
			=
			\left(\sum_{i\in I_E}t_ib_{ij}\right)
			(1-gg_j)=0.
		\end{equation}
		For $j\in I_E$, we have $\chi_j(g)=-1$.
		Multiplying the corresponding equations by $t_j$
		and summing over $j\in I_E$ gives $2v^2=0$.
		Moreover, $hv=\chi(h)vh$ for every $h\in\Gamma$.
		Thus $v$ commutes with every generator of $A$ up to
		multiplication by an element in $k^\times$. Since $v^2=0$, the ideal generated
		by $v$ is square-zero and then $v\in J(A)$. Consequently, $R_E\subseteq J(A)\cap E$.
		
		Conversely, suppose that $0\ne v\in J(A)\cap E$.
		If $N>2$, then $E=kx_i$. Since $J(A)\cap k\Gamma=0$, $a_i(1-g_i^{N_i})=0$, which implies that $a_i=0$.
		Similarly, $c_{ij}(1-g_ig_j)=0$ with $j\ne i$ and hence $c_{ij}=0$.  As a result, $v\in R_E$.

		If $N=2$, write $v=\sum_{i\in I_E}t_ix_i$.
		The left-hand side of \eqref{eq:normal-radical}
		belongs to $J(A)$, whereas its right-hand side
		belongs to $k\Gamma$. For $1\le j\le n$, one has
		\[
		\left(\sum_{i\in I_E}t_ib_{ij}\right)(1-gg_j)=0.
		\]
		That is,  $v\in\ker\Theta_E$.
		Consequently, $J(A)\cap E\subseteq R_E$.
	\end{proof}
	
	\begin{theorem}\label{thm:qls}
		The largest Hopf ideal of $A$ contained in $J(A)$ is $(R)$, that is, $\Imax(A)=(R)$. Furthermore, $B(A^*)=(A/(R))^*$ and
		$A^*$ is generated by its coradical
		exactly when  $R_E=0$ for every nonzero isotypic
		component $E$.
	\end{theorem}
	
	\begin{proof}
		By Corollary~\ref{cor:general-largest} and
		Lemma~\ref{lem:qls-radical-intersections},
		\[
		\Imax(A)
		=\left(\bigoplus_E J(A)\cap E\right)
		=(R).
		\]
		By Lemma~\ref{lem:dual-radical-annihilator}, $B(A^*)=(A/(R))^*$.
		The final assertion follows from
		$B(A^*)=A^*$ if and only if $\Imax(A)=0$.
	\end{proof}

	\begin{proposition}\label{prop:qls-contraction}
		For every nonzero isotypic component
		$F=V_g^\chi$, put $E=E_{g,\chi}$.
		Then
		\[
		\iota(\ker\Psi_F)=R_E=J(A)\cap E.
		\]
		More explicitly, the following hold.
		\begin{itemize}
			\item[(1)] If $F=kx_i$, then
			$\ker\Psi_F=F$ if and only if
			\[
			N_i a_i(1-g_i^{N_i})=0,\qquad
			c_{ij}(1-g_ig_j)=0\quad(j\ne i).
			\]
			
			\item[(2)] If $\dim F>1$, then
			$\operatorname{ord}\chi(g)=2$. For
			$v=\sum_{i\in I_E}t_ix_i$, 
			$v\in\ker\Psi_F$ if and only if
			\[
			\left(\sum_{i\in I_E}t_ib_{ij}\right)
			(1-gg_j)=0,\qquad 1\le j\le n.
			\]
		\end{itemize}
	\end{proposition}
	
	\begin{proof}
		We apply Theorem~\ref{thm:contraction} with
		$W=\mathcal W_{\mathrm{qls}}$ and
		$\kappa=\kappa_D$.
		The space $\mathcal W_{\mathrm{qls}}$ is
		finite-dimensional, homogeneous and generates
		$\mathcal J(V)$. It remains to check its stability
		under $\Aut_{\YD}(V)$.
		
		Every Yetter--Drinfeld automorphism preserves the
		isotypic components. If $F=V_g^\chi$ has dimension
		greater than one, then $\chi(g)=-1$. The quadratic
		relations supported on $F\otimes F$ are spanned by $uv+vu$ with $u,v\in F$.
		For two distinct components
		$F=V_g^\chi$ and $F'=V_h^\psi$, the relations
		between them are spanned by $uw-\psi(g)wu$ with $u\in F$ and $w\in F'$.
		These two relation spaces are preserved by changes
		of bases in the corresponding components. On a
		one-dimensional component $kx_i$, the power
		relation $x_i^{N_i}$ is preserved by rescaling
		$x_i$. Thus $\mathcal W_{\mathrm{qls}}$ is stable
		under $\Aut_{\YD}(V)$. By definition,
		$\kappa_D(\mathcal W_{\mathrm{qls}})\subseteq k\Gamma$.
		
		By Theorem~\ref{thm:contraction}, $\iota(\ker\Psi_F)=J(A)\cap E$.
		Lemma~\ref{lem:qls-radical-intersections} identifies
		$J(A)\cap E$ with $R_E$.
		
		If $F=kx_i$, then $x_i\in\ker\Psi_F$ 
		is described explicitly as follows:
		\[
		N_i a_i(1-g_i^{N_i})=0,\qquad
		c_{ij}(1-g_ig_j)=0\quad(j\ne i).
		\]

		Suppose now that $\dim F>1$. Then
		$\operatorname{ord}\chi(g)=2$. For
		$v=\sum_{i\in I_E}t_ix_i$, the equations
		$\Psi_F(v)=0$ are
		\[
		\left(\sum_{i\in I_E}t_ib_{ij}\right)
		(1-gg_j)=0,\qquad 1\le j\le n.
		\]
		By the definition of $\Theta_E$, these are exactly
		the equations $\Theta_E(v)=0$. Hence
		$\iota(\ker\Psi_F)=\ker\Theta_E=R_E$.
	\end{proof}
	
	\begin{corollary}\label{cor:qls-parameter-criterion}
		The algebra $A^*$ is generated by its coradical
		if and only if the following conditions hold:
		\begin{itemize}
			\item[(1)] For every one-dimensional isotypic
			component $E_{g_i,\chi_i}=kx_i$, either
			$a_i\ne0$ or $c_{ij}\ne0$ for some $j\ne i$.
			
			\item[(2)] For every isotypic component
			$E=E_{g,\chi}$ of dimension greater than one,
			the matrix in \eqref{eq:pairing-map}, after
			deleting the rows indexed by $j$ with
			$gg_j=1$, has full column rank.
		\end{itemize}
	\end{corollary}
	
	\begin{proof}
		By Theorem~\ref{thm:qls}, $A^*$ is generated by
		its coradical if and only if $R_E=0$ for every
		nonzero isotypic component $E$.
		For one-dimensional components, Proposition
		\ref{prop:qls-contraction}(1) gives condition (1).
		For components of dimension greater than one,
		Proposition~\ref{prop:qls-contraction}(2) identifies
		$R_E$ with the kernel of the matrix in
		\eqref{eq:pairing-map}. Its vanishing is equivalent
		to the full column rank condition in (2).
	\end{proof}
	
	\begin{remark}\label{rmk:example-sec6}
		The nonvanishing of all power and linking parameters
		does not necessarily imply that
		$J(A)\cap E_{g,\chi}=0$ when $\dim E_{g,\chi}>1$.
		
		Let $\Gamma=C_4=\langle g\rangle$ and choose
		$\chi\in\widehat{\Gamma}$ with $\chi(g)=-1$.
		Take $F=V_g^\chi$ with basis $u,v$, and choose
		the compatible datum $D$ given by
		$a_1=a_2=1$ and $c_{12}=2$.
		Let $A=A_D$, and write $x=\iota(u)$ and
		$y=\iota(v)$. The corresponding relations in $A$ are
		\[
		x^2=y^2=1-g^2,\qquad
		xy+yx=2(1-g^2).
		\]
		With respect to the basis $\{u,v\}$, the matrix
		in \eqref{eq:pairing-map} is
		\[
		\begin{pmatrix}
			2&2\\
			2&2
		\end{pmatrix}.
		\]
		Thus $\ker\Psi_F=k(u-v)$. By Proposition~\ref{prop:qls-contraction}, $J(A)\cap E_{g,\chi}=\iota(\ker\Psi_F)=k(x-y)$.
	\end{remark}

	\subsubsection{Isomorphisms and changes of generators}
	
	Let $A$ and $A'$ be liftings of quantum linear spaces
	with data $(\Gamma,g_i,\chi_i)$ and
	$(\Gamma',g'_j,\chi'_j)$, respectively.
	
	\begin{proposition}\label{prop:qls-isomorphism}
		
		Let $\varphi:A\to A'$ be a Hopf algebra isomorphism. Then $\beta=\varphi\mid_{\Gamma}$ is  a group isomorphism
		$\beta:\Gamma\to\Gamma'$ and $\varphi$ satisfies
		\begin{equation}\label{eq:datum-iso}
			\varphi(x_i)=\sum_j M_{ji}x'_j,\qquad
			M_{ji}=0\quad\text{unless}\quad
			\beta(g_i)=g'_j,\quad
			\chi_i=\chi'_j\circ\beta.
		\end{equation}
		Conversely, any algebra map defined by
		\eqref{eq:datum-iso} with $M$ invertible is a
		Hopf algebra isomorphism.
	\end{proposition}
	
	\begin{proof}
		Clearly, $\beta=\varphi|_\Gamma$ is a group
		isomorphism. For each $i$, we have
		$\varphi(x_i)\in\Pp_{\beta(g_i),1}(A')$, and for $h\in\Gamma$,
		\[
		\beta(h)\varphi(x_i)\beta(h)^{-1}
		=\chi_i(h)\varphi(x_i).
		\]
		By \cite[Lemma 5.4]{ASLifting}, the quotient
		$\Pp_{\beta(g_i),1}(A')/k(1-\beta(g_i))$
		is spanned by the classes of those $x'_j$ for which
		$g'_j=\beta(g_i)$. Since $\chi_i\ne\varepsilon$,
		the component of $\varphi(x_i)$ in
		$k(1-\beta(g_i))$ is zero. The conjugation relation
		shows that only those $x'_j$ satisfying
		$\chi_i=\chi'_j\circ\beta$ can occur and so
		\eqref{eq:datum-iso} follows.
		
		Conversely, the conditions in \eqref{eq:datum-iso}
		ensure that $\varphi$ is a coalgebra map and then is a Hopf algebra map. Since $\beta$
		and $M$ are invertible, $\varphi$ is surjective.
		The matching of the pairs $(g_i,\chi_i)$
		preserves the corresponding orders $N_i$, so $\dim A=\dim A'$. Therefore, 
		$\varphi$ is an isomorphism.
	\end{proof}
	
	We restrict to $n=2$ and write $x=x_1$, $y=x_2$,
	$a=a_1$, $b=a_2$ and $c=c_{12}$.
	Fix the ordered data $(g_i,\chi_i)$ and consider
	isomorphisms inducing $\beta=\id_\Gamma$.
	By Proposition~\ref{prop:qls-isomorphism}, if
	$(g_1,\chi_1)\ne(g_2,\chi_2)$, the only possible
	changes of generators are $x\mapsto sx$ and
	$y\mapsto ty$, with $s,t\in k^\times$. Accordingly,
	the parameters transform as
	\[
	(a,b,c)\longmapsto
	(s^{N_1}a,t^{N_2}b,stc).
	\]
	By Corollary~\ref{cor:qls-parameter-criterion}, $A^*$
	is generated by its coradical if and only if $c\ne0$
	or $ab\ne0$. This condition is preserved by these
	rescalings.
	
	Suppose now that
	$(g_1,\chi_1)=(g_2,\chi_2)=(g,\chi)$.
	Then $N_1=N_2=2$, and every change of basis in
	$E=\operatorname{span}\{x,y\}$ preserves the structure of 
	Yetter--Drinfeld modules. Associate to the parameters
	$(a,b,c)$ the symmetric matrix
	\[
	B(a,b,c)=
	\begin{pmatrix}
		2a&c\\
		c&2b
	\end{pmatrix}.
	\]
	For $v=t_1x+t_2y$, the map in \eqref{eq:pairing-map}
	is
	\[
	\Theta_E(v)=
	\bigl((2at_1+ct_2)(1-g^2),\,
	(ct_1+2bt_2)(1-g^2)\bigr).
	\]
	If $g^2\ne1$, then $1-g^2\ne0$ in $k\Gamma$,
	and hence $\ker\Theta_E$ is identified with $\ker B(a,b,c)$
	in the basis $\{x,y\}$.
	If $g^2=1$, compatibility forces $a=b=c=0$, so
	both maps vanish. Since $E$ is the only nonzero isotypic component,
	Proposition~\ref{prop:qls-contraction} and
	Theorem~\ref{thm:qls} give $\dim(J(A)\cap E)=2-\rank B(a,b,c)$ and moreover, 
	$B(A^*)=A^*$ if and only if $\rank B(a,b,c)=2$.

	\begin{proposition}\label{prop:qls-congruence}
		Two parameter triples $(a,b,c)$ and
		$(a',b',c')$ define isomorphic Hopf algebras
		(with the same ordered data) if and only if
		$B(a,b,c)$ and $B(a',b',c')$ have the same rank.
	\end{proposition}
	
	\begin{proof}
		If $g^2=1$, compatibility forces both triples to be zero.
		Suppose that $g^2\ne1$. For $M\in\GL_2(k)$, let
		$x'_a=\sum_{i=1}^2M_{ia}x_i$. From the defining relations, one has
		\[
		x'_ax'_b+x'_bx'_a
		=(M^{\mathsf t}B(a,b,c)M)_{ab}(1-g^2).
		\]
		By Proposition~\ref{prop:qls-isomorphism}, every Hopf
		algebra isomorphism between the two liftings fixes $g$
		and has this form on the generators. After the change of basis, $B(a,b,c)$ is replaced
		by $M^{\mathsf t}B(a,b,c)M$. Conversely, every
		such congruence is induced by a change of basis.
		Observe that $k$ is algebraically closed of characteristic
		zero. Two symmetric matrices of size two are
		congruent if and only if they have the same rank.
	\end{proof}
	
	\begin{remark}
		For the example listed in Remark \ref{rmk:example-sec6}, the change of generators
		$x'=x-y$, $y'=y$ replaces $(1,1,2)$ by $(0,1,0)$.
		Both parameter matrices have rank one, although the
		nonvanishing of individual parameters changes.
	\end{remark}

	\subsection{Liftings of finite Cartan type}\label{sec:cartan}
	
	\begin{definition}\cite[Definition 2.1]{ASCartan}\label{def:cartan-datum}
		A datum of Cartan type
		\[
		\mathcal D=\mathcal D(\Gamma,(g_i)_{1\le i\le\theta},
		(\chi_i)_{1\le i\le\theta},(a_{ij})_{1\le i,j\le\theta})
		\]
		consists of an abelian group $\Gamma$, elements
		$g_i\in\Gamma$, characters $\chi_i\in\widehat\Gamma$,
		$1\le i\le\theta$, and a generalized Cartan matrix
		$(a_{ij})$ of size $\theta$ satisfying
		\[
		q_{ij}q_{ji}=q_{ii}^{a_{ij}},\quad q_{ii}\ne1,
		\qquad\text{with }q_{ij}=\chi_j(g_i)
		\quad(1\le i,j\le\theta).
		\]
		It is of finite Cartan type if $(a_{ij})$ is of finite type.
	\end{definition}
	
	Fix a datum $\mathcal D$ of finite Cartan type with $\Gamma$
	finite. Put $I=\{1,2,\ldots,\theta\}$, and let
	$\mathcal X=\{I_1,\ldots,I_t\}$ be the set of connected
	components of its Dynkin diagram. Write $i\sim j$ if $i$ and $j$
	belong to the same component. We assume the order conditions
	\cite[(3.1)--(3.2)]{ASCartan}:
	\[
	\begin{gathered}
		q_{ii}\text{ has odd order},\\
		\text{the order of }q_{ii}\text{ is prime to }3
		\text{ if }i\text{ lies in a component }G_2.
	\end{gathered}
	\]
	For $J\in\mathcal X$, let $N_J$ be the common order of
	$q_{ii}$, $i\in J$. For the application below, the present paper
	imposes the additional assumption
	\[
	(g_i,\chi_i)\ne(g_j,\chi_j)\qquad(i\ne j).
	\]
	
	\begin{remark}
		Quantum linear spaces are  the Nichols algebras
		of Cartan type
		$A_1\times\cdots\times A_1$.
		The preceding subsection treats their liftings without
		imposing the order conditions or the distinctness
		assumption used here. It also provides additional
		structural information that will be used in the next section.
	\end{remark}
	
	\begin{definition}\cite[Section 2.1]{ASCartan}\label{def:cartan-roots}
		Let $\mathbb Z[I]$ be the root lattice with basis
		$\alpha_1,\ldots,\alpha_\theta$.
		The Weyl group $W\subseteq\Aut(\mathbb Z[I])$ is generated
		by the reflections
		$s_i(\alpha_j)=\alpha_j-a_{ij}\alpha_i$.
		The root system and its positive roots are
		\[
		\Phi=\bigcup_{i=1}^{\theta}W(\alpha_i),\qquad
		\Phi^+=\Bigl\{\alpha\in\Phi:
		\alpha=\sum_{i=1}^{\theta}n_i\alpha_i,\quad
		n_i\ge0\text{ for all }i\Bigr\}.
		\]
		For $J\in\mathcal X$, let $\Phi_J$ be the root system of
		$(a_{ij})_{i,j\in J}$, so that
		$\Phi^+=\bigsqcup_{J\in\mathcal X}\Phi_J^+$.
		For $\alpha=\sum_{i=1}^{\theta}n_i\alpha_i\in\mathbb Z[I]$,
		set
		\[
		g_\alpha=g_1^{n_1}g_2^{n_2}\cdots g_\theta^{n_\theta},
		\qquad
		\chi_\alpha=\chi_1^{n_1}\chi_2^{n_2}\cdots\chi_\theta^{n_\theta},
		\]
		and $q_{\alpha\beta}=\chi_\beta(g_\alpha)$ for
		$\alpha,\beta\in\mathbb Z[I]$. For $\alpha\in\Phi^+$,
		we also write $\supp\alpha=\{i:n_i>0\}$.
	\end{definition}
	
	We first recall the construction of \cite[Sections 2.1--2.2]{ASCartan}
	for a connected datum $\mathcal D$. Until
	Definition~\ref{def:cartan-root-corrections}, all notation refers
	to this connected datum. In particular,
	\[
	N=\text{order of }q_{ii},\quad1\le i\le\theta,
	\]
	as in \cite[(2.6)]{ASCartan}. Fix a reduced decomposition
	of the longest element of $W$,
	\[
	w_0=s_{i_1}s_{i_2}\cdots s_{i_p}.
	\]
	The corresponding convex ordering of $\Phi^+$ is
	\[
	\beta_l=s_{i_1}\cdots s_{i_{l-1}}(\alpha_{i_l}),
	\qquad1\le l\le p.
	\]
	
	\begin{definition}\cite[Definition 2.4]{ASCartan}\label{def:cartan-R}
		Let $V=V(\mathcal D)\in\YD$
		have basis $x_1,\ldots,x_\theta$, with
		$x_i\in V_{g_i}^{\chi_i}$. Define
		\[
		R(\mathcal D)=T(V)/
		\bigl((\ad_cx_i)^{1-a_{ij}}(x_j)
		\mid1\le i\ne j\le\theta\bigr).
		\]
	\end{definition}
	
	For $\alpha\in\Phi^+$, let $x_\alpha\in R(\mathcal D)$
	be the root vector defined in \cite[Section 2.1]{ASCartan}
	by iterated braided commutators for the chosen convex ordering.
	We use the same notation for the corresponding polynomial in
	$k\langle x_1,\ldots,x_\theta\rangle$ and its images in the
	quotients below. If $\alpha=\sum_i n_i\alpha_i$, then
	$x_\alpha$ is multihomogeneous in the simple generators of
	multidegree $(n_1,\ldots,n_\theta)$.
	
	\begin{definition}\cite[Definitions 2.5 and 2.7]{ASCartan}\label{def:cartan-power-data}
		Let $K(\mathcal D)$ be the subalgebra of $R(\mathcal D)$
		generated by $x_\alpha^N$, $\alpha\in\Phi^+$.
		For $1\le l\le p$ and $a=(a_1,\ldots,a_p)\in\mathbb N^p$,
		where $\mathbb N=\{0,1,2,\ldots\}$, set
		\[
		h_l=g_{\beta_l}^N,\qquad
		\eta_l=\chi_{\beta_l}^N,\qquad
		z_l=x_{\beta_l}^N,
		\]
		\[
		\begin{aligned}
			z^a&=z_1^{a_1}z_2^{a_2}\cdots z_p^{a_p}\in K(\mathcal D),\\
			h^a&=h_1^{a_1}h_2^{a_2}\cdots h_p^{a_p}\in\Gamma,\\
			\eta^a&=\eta_1^{a_1}\eta_2^{a_2}\cdots\eta_p^{a_p}
			\in\widehat\Gamma,\\
			\underline a&=a_1\beta_1+a_2\beta_2+\cdots+a_p\beta_p
			\in\mathbb Z[I].
		\end{aligned}
		\]
		For $\alpha=\sum_{i=1}^{\theta}n_i\alpha_i\in\mathbb Z[I]$,
		let $\operatorname{ht}(\alpha)=\sum_{i=1}^{\theta}n_i$,
		and let $e_l=(\delta_{kl})_{1\le k\le p}\in\mathbb N^p$.
	\end{definition}
	
	The element $z^a$ has $\mathbb Z[I]$-degree $N\underline a$
	and belongs to $K(\mathcal D)_{h^a}^{\eta^a}$.
	If $\underline b+\underline c=\underline a$, then
	$h^a=h^bh^c$ and $\eta^a=\eta^b\eta^c$; if also $b,c\ne0$,
	both $\underline b$ and $\underline c$ have smaller height
	than $\underline a$.
	
	\begin{proposition}\cite[Theorem 2.6, Lemma 2.8 and (2.17)]{ASCartan}\label{prop:cartan-power-coalgebra}
		The algebra $K(\mathcal D)$ is a braided Hopf subalgebra
		of $R(\mathcal D)$, and the elements $z^ag$,
		$a\in\mathbb N^p$, $g\in\Gamma$, form a basis of
		$K(\mathcal D)\#k\Gamma$. Its multiplication satisfies
		\[
		z^bz^c=\gamma_{b,c}z^{b+c},\qquad
		\gamma_{b,c}=\prod_{k>l}\eta_l(h_k)^{b_kc_l},
		\qquad
		h^az^b=\eta^b(h^a)z^bh^a.
		\]
		For $0\ne a\in\mathbb N^p$, the elements $t_{b,c}^a$ in $k$
		are uniquely determined by
		\begin{equation}\label{eq:cartan-braided-coproduct}
			\Delta_{K(\mathcal D)}(z^a)
			=z^a\otimes1+1\otimes z^a+
			\sum_{b,c\ne0,\ \underline b+\underline c=\underline a}
			t_{b,c}^a\,z^b\otimes z^c.
		\end{equation}
		In the bosonization,
		\begin{equation}\label{eq:cartan-power-coproduct}
			\Delta_{K(\mathcal D)\#k\Gamma}(z^a)
			=h^a\otimes z^a+z^a\otimes1+
			\sum_{b,c\ne0,\ \underline b+\underline c=\underline a}
			t_{b,c}^a\,z^bh^c\otimes z^c.
		\end{equation}
		Here and below, all multiindices in these sums belong to
		$\mathbb N^p$.
	\end{proposition}
	
	For $n\ge0$, let $K(\mathcal D)_n$ be the span of the $z^a$
	with $\operatorname{ht}(\underline a)\le n$.
	Then $K(\mathcal D)_n\#k\Gamma$ is a subcoalgebra of
	$K(\mathcal D)\#k\Gamma$.
	
	\begin{lemma}\cite[Lemmas 2.9--2.10]{ASCartan}\label{lem:cartan-normalized-recursion}
		Let $n\ge1$, and let
		$(\mu_a)_{0\ne a\in\mathbb N^p,\ \operatorname{ht}(\underline a)\le n}$
		be a family of elements in $k$ satisfying $\mu_a=0$ if $h^a=1$.
		Define $u^a\in k\Gamma$ by induction on
		$\operatorname{ht}(\underline a)$ through
		\begin{equation}\label{eq:cartan-normalized-recursion}
			u^a=\mu_a(1-h^a)+
			\sum_{b,c\ne0,\ \underline b+\underline c=\underline a}
			t_{b,c}^a\,\mu_bu^c.
		\end{equation}
		The map
		\[
		\varphi:K(\mathcal D)_n\#k\Gamma\longrightarrow k\Gamma,
		\qquad \varphi|\Gamma=\id,\qquad
		\varphi(z^ag)=u^ag
		\]
		for $0\ne a\in\mathbb N^p$ with
		$\operatorname{ht}(\underline a)\le n$ and $g\in\Gamma$
		is a coalgebra map. In particular,
		\begin{equation}\label{eq:cartan-u-coproduct}
			\Delta(u^a)=h^a\otimes u^a+u^a\otimes1+
			\sum_{b,c\ne0,\ \underline b+\underline c=\underline a}
			t_{b,c}^a\,u^bh^c\otimes u^c.
		\end{equation}
		The family $(\mu_a)$ associated to this map is unique
		under the condition $\mu_a=0$ if $h^a=1$.
	\end{lemma}
	
	\begin{definition}\cite[Section 2.2, preceding Theorem 2.13]{ASCartan}\label{def:cartan-M}
		Let $M(\mathcal D)$ be the set of families
		$\mu=(\mu_l)_{1\le l\le p}$ of elements in $k$ such that
		\[
		\text{if }h_l=1\text{ or }\eta_l\ne\eps,
		\text{ then }\mu_l=0.
		\]
	\end{definition}
	
	\begin{theorem}\cite[Lemma 2.12(2) and Theorem 2.13]{ASCartan}\label{thm:cartan-parameter-map}
		For $\mu=(\mu_l)_{1\le l\le p}\in M(\mathcal D)$,
		there is exactly one Hopf algebra map
		\[
		\varphi_\mu:K(\mathcal D)\#k\Gamma\longrightarrow k\Gamma,
		\qquad\varphi_\mu|\Gamma=\id,
		\]
		whose associated family $(\mu_a)_{0\ne a\in\mathbb N^p}$
		in Lemma~\ref{lem:cartan-normalized-recursion} satisfies
		$\mu_{e_l}=\mu_l$ for $1\le l\le p$.
		The assignment $\mu\mapsto\varphi_\mu$ is a bijection
		from $M(\mathcal D)$ to the Hopf algebra maps that are the
		identity on $\Gamma$. For $u^a=\varphi_\mu(z^a)$ and
		$u_l=u^{e_l}$, one has
		\begin{equation}\label{eq:coproduct_formu}
			u^a=\prod_{a_l>0}u_l^{a_l},\qquad
			\eta^a\ne\eps\ \Longrightarrow\ \mu_a=0,\ u^a=0
			\quad(0\ne a\in\mathbb N^p).
		\end{equation}
	\end{theorem}
	
	The construction in the proof of \cite[Theorem 2.13]{ASCartan}
	is inductive in $\operatorname{ht}(\underline a)$. If $a=e_l$,
	it uses \cite[(2.19)]{ASCartan}:
	\begin{equation}\label{eq:cartan-generator-recursion}
		u^a=\mu_l(1-h^a)+
		\sum_{b,c\ne0,\ \underline b+\underline c=\underline a}
		t_{b,c}^a\mu_bu^c.
	\end{equation}
	If $a=(a_1,\ldots,a_l,0,\ldots,0)$, $a_l\ge1$, and
	$a\ne e_l$, write $a=r+s$, where $0\ne r$ and $s=e_l$,
	and set $u^a=u^ru^s$. The associated elements $\mu_a$ in $k$ are
	then determined by \eqref{eq:cartan-normalized-recursion}
	and the condition $\mu_a=0$ if $h^a=1$.
	The parameters $\mu_l$ are the prescribed values $\mu_{e_l}$;
	the other $\mu_a$ are determined by this construction.
	
	\begin{definition}\cite[Definition 2.14]{ASCartan}\label{def:cartan-root-corrections}
		For $\mu\in M(\mathcal D)$, define
		\[
		u_l(\mu)=\varphi_\mu(z_l)\in k\Gamma,\qquad1\le l\le p.
		\]
		If $\alpha=\beta_l\in\Phi^+$, set
		$u_\alpha(\mu)=u_l(\mu)$.
	\end{definition}
	
	Each $u_\alpha(\mu)$ lies in the augmentation ideal of
	$k[g_i^N\mid1\le i\le\theta]$ by
	\cite[(2.15)]{ASCartan}.
	
	We now return to an arbitrary datum $\mathcal D$ with
	components $\mathcal X=\{I_1,\ldots,I_t\}$.
	For $J\in\mathcal X$, denote its restriction to $J$ by
	$\mathcal D_J$, as in the proof of \cite[Theorem 4.5]{ASCartan}.
	The preceding connected construction is applied to
	$\mathcal D_J$ with $N=N_J$, retaining its local notation.
	Choose a reduced decomposition of the longest element
	$w_{0,J}$ of the Weyl group $W_J$ of $\Phi_J$. Let
	$p_J=|\Phi_J^+|$ and let
	\[
	\Phi_J^+=\{\beta_{J,1},\ldots,\beta_{J,p_J}\}
	\]
	be the corresponding convex ordering. As in
	\cite[Section 3.1]{ASCartan}, the reduced decomposition
	$w_0=w_{0,I_1}w_{0,I_2}\cdots w_{0,I_t}$ gives the ordering
	\[
	\Phi^+=\{\beta_{I_1,1},\ldots,\beta_{I_1,p_{I_1}},\ldots,\beta_{I_t,1},\ldots,\beta_{I_t,p_{I_t}}\}=\{\beta_1,\ldots,\beta_p\},\quad
	p=\sum_{J\in\mathcal X}p_J.
	\]
	Let $V=V(\mathcal D)\in\YD$
	have basis $x_1,\ldots,x_\theta$, with $x_i\in V_{g_i}^{\chi_i}$.
	
	\begin{definition}\cite[Definition 3.1 and (3.4)]{ASCartan}\label{def:cartan-parameters}
		A family $\lambda=(\lambda_{ij})_{1\le i<j\le\theta,\ i\not\sim j}$
		of elements in $k$ is a family of linking parameters for
		$\mathcal D$ if
		\[
		\text{if }g_ig_j=1\text{ or }\chi_i\chi_j\ne\eps,
		\text{ then }\lambda_{ij}=0.
		\]
		Extend the family by
		\[
		\lambda_{ji}=-q_{ji}\lambda_{ij}
		\qquad(1\le i<j\le\theta,\ i\not\sim j).
		\]
	\end{definition}
	
	\begin{definition}\cite[Definition 3.2]{ASCartan}\label{def:cartan-U}
		Let $\lambda$ be a family of linking parameters for
		$\mathcal D$. The Hopf algebra $U(\mathcal D,\lambda)$
		is the quotient of
		$k\langle x_1,\ldots,x_\theta\rangle\#k\Gamma$
		by the ideal generated by the Serre and linking relations
		\[
		\begin{aligned}
			&\ad_c(x_i)^{1-a_{ij}}(x_j),
			&&1\le i,j\le\theta,\ i\sim j,\ i\ne j,\\
			&x_ix_j-q_{ij}x_jx_i-\lambda_{ij}(1-g_ig_j),
			&&1\le i<j\le\theta,\ i\not\sim j.
		\end{aligned}
		\]
		The images of $x_i$ and $g\in\Gamma$ are denoted by the
		same symbols, and
		\[
		\Delta(x_i)=g_i\otimes x_i+x_i\otimes1,
		\qquad\Delta(g)=g\otimes g.
		\]
	\end{definition}
	
	\begin{proposition}\cite[Theorem 3.3(1)]{ASCartan}\label{prop:cartan-U-pbw}
		The elements
		\[
		x_{\beta_1}^{a_1}x_{\beta_2}^{a_2}\cdots x_{\beta_p}^{a_p}g,
		\qquad a_1,a_2,\ldots,a_p\ge0,\quad g\in\Gamma,
		\]
		form a basis of $U(\mathcal D,\lambda)$.
	\end{proposition}
	
	\begin{definition}\cite[Definition 4.3]{ASCartan}\label{def:cartan-lifting}
		A family $\mu=(\mu_\alpha)_{\alpha\in\Phi^+}$ of elements
		in $k$ is a family of root vector parameters for $\mathcal D$
		if, for $\alpha\in\Phi_J^+$ and $J\in\mathcal X$,
		\[
		\text{if }g_\alpha^{N_J}=1\text{ or }\chi_\alpha^{N_J}\ne\eps,
		\text{ then }\mu_\alpha=0.
		\]
		For $J\in\mathcal X$ and $\alpha\in\Phi_J^+$, define
		\[
		\pi_J(\mu)=(\mu_\beta)_{\beta\in\Phi_J^+},\qquad
		u_\alpha(\mu)=u_\alpha(\pi_J(\mu)).
		\]
		Here Definition~\ref{def:cartan-root-corrections} is applied
		to $\mathcal D_J$: its local $\beta_l$ is $\beta_{J,l}$,
		and its prescribed parameter $\mu_l$ is $\mu_{\beta_{J,l}}$.
		For linking parameters $\lambda$, set
		\begin{equation}\label{eq:cartan-powers}
			u(\mathcal D,\lambda,\mu)=U(\mathcal D,\lambda)/
			\bigl(x_\alpha^{N_J}-u_\alpha(\mu)
			\mid\alpha\in\Phi_J^+,\ J\in\mathcal X\bigr).
		\end{equation}
	\end{definition}
	
	Write $A=u(\mathcal D,\lambda,\mu)$, and retain the symbols
	$x_i$ and $g$ for their images in $A$. For $1\le l\le p$,
	let $N_l=N_J$ if $\beta_l\in\Phi_J^+$.
	
	\begin{theorem}\cite[Theorem 4.5]{ASCartan}\label{thm:cartan-pbw}\label{prop:cartan-pbw}
		The algebra $u(\mathcal D,\lambda,\mu)$ is a finite-dimensional
		pointed quotient Hopf
		algebra of $U(\mathcal D,\lambda)$ with
		$G(u(\mathcal D,\lambda,\mu))\cong\Gamma$. The elements
		\[
		x_{\beta_1}^{a_1}x_{\beta_2}^{a_2}\cdots x_{\beta_p}^{a_p}g,
		\qquad0\le a_l<N_l,\quad1\le l\le p,\quad g\in\Gamma,
		\]
		form a basis of $u(\mathcal D,\lambda,\mu)$. In particular,
		\[
		\dim u(\mathcal D,\lambda,\mu)
		=\prod_{J\in\mathcal X}N_J^{|\Phi_J^+|}|\Gamma|.
		\]
	\end{theorem}
	
	We identify $k\Gamma$ with its image in $A$.
	
	\begin{theorem}\cite[Theorem 5.1]{ASCartan}\label{thm:cartan-nichols}
		The Nichols algebra $\B(V)$ is the quotient of
		$T(V)$ by the ideal generated by
		\[
		\begin{aligned}
			&\ad_c(x_i)^{1-a_{ij}}(x_j),
			&&1\le i,j\le\theta,\ i\ne j,\\
			&x_\alpha^{N_J},
			&&\alpha\in\Phi_J^+,\ J\in\mathcal X.
		\end{aligned}
		\]
	\end{theorem}
	
	\begin{corollary}\cite[Corollary 5.2]{ASCartan}\label{cor:cartan-gr}\label{prop:cartan-gr}
		$\gr(u(\mathcal D,\lambda,\mu))\cong u(\mathcal D,0,0)\cong\B(V)\#k\Gamma$.
	\end{corollary}
	
	\begin{lemma}\label{lem:cartan-root-parameter}
		For every family of root vector parameters $\mu$,
		\[
		\bigcup_{\beta:\,u_\beta(\mu)\ne0}\supp\beta
		=\bigcup_{\beta:\,\mu_\beta\ne0}\supp\beta.
		\]
	\end{lemma}
	
	\begin{proof}
		It suffices to work in each connected component.
		Use the connected notation of \cite[Section 2.2]{ASCartan},
		so that $\Phi^+=\{\beta_1,\ldots,\beta_p\}$,
		$\mu_l=\mu_{\beta_l}$, $u^a=\varphi_\mu(z^a)$,
		and $u_l=u^{e_l}=u_{\beta_l}(\mu)$.
		Fix a simple index $i$.
		
		Suppose that $\mu_l=0$ whenever
		$i\in\supp\beta_l$. We prove by induction on
		$\operatorname{ht}(\underline a)$ that $u^a=\mu_a=0$
		if the coefficient of $\alpha_i$ in $\underline a$ is
		positive. In every summand of
		\eqref{eq:cartan-normalized-recursion}, this coefficient
		is positive in $\underline b$ or $\underline c$.
		Both have smaller height and then by induction, $\mu_bu^c=0$.
		If $a=e_l$, then by hypothesis $\mu_a=\mu_l=0$,
		and by recursion, $u^a=0$.
		If $a\ne e_l$ for every $l$, choose $l$ with $a_l>0$
		and $i\in\supp\beta_l$. Then
		$\operatorname{ht}(\beta_l)<\operatorname{ht}(\underline a)$,
		and by induction, $u_l=0$.  By \eqref{eq:coproduct_formu}, $u^a=0$.
		The recursion now gives $\mu_a(1-h^a)=0$.
		If $h^a\ne1$, then $\mu_a=0$; otherwise, it follows from the normalization.
		
		Conversely, suppose that $u_l=0$ whenever
		$i\in\supp\beta_l$. Then \eqref{eq:coproduct_formu} implies that $u^a=0$
		whenever the coefficient of $\alpha_i$ in $\underline a$
		is positive. By induction in
		\eqref{eq:cartan-normalized-recursion}, $\mu_a(1-h^a)=0$, since each summand has $\mu_b=0$
		by induction or $u^c=0$ by \eqref{eq:coproduct_formu}.
		The normalization again gives $\mu_a=0$.
		Taking $a=e_l$ proves the assertion in each component.
	\end{proof}

	\begin{theorem}\label{thm:cartan}
		$A^*$ is generated by the coradical
		if and only if
		\[
		\bigcup_{\substack{\beta\in\Phi_J^+,\ J\in\mathcal X\\
				u_\beta(\mu)\ne0}}
		\supp\beta
		\ \cup\
		\bigcup_{\substack{i<j,\ i\not\sim j\\
				\lambda_{ij}(1-g_ig_j)\ne0}}
		\{i,j\}
		=
		\{1,\ldots,\theta\}.
		\]
	\end{theorem}
	
	\begin{proof}
		Suppose that the displayed union equals $\{1,\ldots,\theta\}$.
		The pairs $(g_i,\chi_i)$ are distinct. Then every nonzero space $E_{g,\chi}$ is
		one-dimensional. By Corollary~\ref{cor:equal-type},
		if $A^*$ is not generated by its coradical,
		then there is some $i$ such that $x_i\in J(A)$.
		
		Suppose that $i\in\supp\beta$ and
		$u_\beta(\mu)\ne0$. Since every monomial in $x_\beta$
		contains $x_i$, it follows that $x_\beta\in J(A)$.
		From the relation $x_\beta^{N_J}=u_\beta(\mu)$, we have  $u_\beta(\mu)\in J(A)\cap k\Gamma=0$,
		a contradiction.
		
		Likewise, if $i\in\{r,s\}$ and
		$\lambda_{rs}(1-g_rg_s)\ne0$, then 
		\[
		x_rx_s-q_{rs}x_sx_r=\lambda_{rs}(1-g_rg_s)\in J(A)\cap k\Gamma=0,
		\]
		which is a contradiction.
		Consequently, $A^*$ is generated by its coradical.
		
		Conversely, suppose that $i$ does not belong to
		the displayed union. Assign degree one to $x_i$
		and degree zero to all other generators and to
		$\Gamma$. The Serre relations are homogeneous for
		this degree.
		
		For $\beta=\sum_jm_j\alpha_j\in\Phi_J^+$,
		the element $x_\beta^{N_J}$ has degree $N_Jm_i$.
		If $m_i>0$, the assumption on $i$ implies
		$u_\beta(\mu)=0$. If $m_i=0$, both sides of the
		root vector relation have degree zero.
		
		For a linking relation involving $x_i$, the
		assumption on $i$ forces
		$\lambda_{rs}(1-g_rg_s)=0$. Every other linking
		relation has degree zero.
		
		Consequently, all defining relations of $A$
		are homogeneous, and the presentation defines
		a nonnegative grading of $A$. The coproduct and
		antipode preserve the grading, since
		\[
		\Delta(x_j)=x_j\otimes1+g_j\otimes x_j,
		\qquad
		S(x_j)=-g_j^{-1}x_j,
		\]
		and all elements of $\Gamma$ have degree zero.
		Thus $A$ is a graded Hopf algebra.
		
		Let $A_+$ be the ideal of positive-degree elements.
		Since $A$ is finite-dimensional, $A_+$ is nilpotent.
		In particular, $x_i\in J(A)$. By
		Lemma~\ref{lem:pointed-detect}, $A^*$ is not
		generated by its coradical.
	\end{proof}

	Hu and Pei studied unified presentations and cocycle
	deformations of two-parameter quantum groups
	\cite{HuPeiI,HuPeiII}. Quantum-shuffle and cotensor
	constructions of multiparameter quantum groups appear
	in \cite{PeiHuRosso,LiHuRosso}, while convex PBW bases
	and restricted two-parameter quantum groups of types
	$B$, $G_2$ and $F_4$ were constructed in
	\cite{HuWangB,HuWangG2,ChenHuWangF4}.
	We consider the restricted multiparameter quantum
	groups of García--Gavarini.
	
	\begin{corollary}\label{cor:restricted-multiparameter}
		Let $\mathfrak g$ be a semisimple Lie algebra of rank $n$,
		with Cartan matrix $(a_{ij})_{1\le i,j\le n}$ and
		symmetrizing integers $d_i$. Let $q\in k$ be a primitive
		$\ell$-th root of unity, where $\ell>1$ is odd and
		coprime to every $d_i$. Set
		$\mathbf q=(q_{ij})=(q^{b_{ij}})$, where $b_{ij}\in\mathbb Z$,
		$b_{ii}=2d_i$ and $b_{ij}+b_{ji}=2d_i a_{ij}$.
		If $A$ is the corresponding restricted small multiparameter
		quantum group over $k$ of \cite[Section 7.3]{GarciaGavarini},
		then $A^*$ is generated by the coradical.
	\end{corollary}
	
	\begin{proof}
		From \cite[Definition 3.1.1, Theorems 7.3.1 and 7.3.3]
		{GarciaGavarini}, $A$ is a lifting of a Nichols algebra
		of finite Cartan type with $\theta=2n$ and $\mu=0$.
		The divided powers and toral quantum binomial coefficients
		of degree less than $\ell$ lie in the algebra generated by
		$E_i,F_i,K_i^{\pm1},L_i^{\pm1}$, since their defining
		denominators are nonzero in $k$. More precisely, for $1\le i\le n$, take
		\[
		x_i=E_i,\quad x_{n+i}=F_iL_i^{-1},\quad
		g_i=K_i,\quad g_{n+i}=L_i^{-1}.
		\]
		where $\Delta(x_r)=x_r\otimes1+g_r\otimes x_r$ for $1\le r\le 2n$.
		With respect to the ordered generators
		$x_1,\ldots,x_{2n}$, the braiding matrix is
		\[
		\widetilde{\mathbf q}=
		\begin{pmatrix}
			\mathbf q & \mathbf q^{[-1]}\\
			\mathbf q^{\mathsf t} & (\mathbf q^{[-1]})^{\mathsf t}
		\end{pmatrix},
		\qquad
		\mathbf q^{[-1]}=(q_{ij}^{-1})_{1\le i,j\le n}.
		\]
		Thus the Cartan matrix consists of two copies of
		$(a_{ij})$, with diagonal braiding entries
		$q_{ii}$ and $q_{ii}^{-1}$, both of order $\ell$.
		The PBW basis shows that
		\[
		\Gamma=\langle K_1,\ldots,K_n,L_1,\ldots,L_n\rangle
		\simeq(\mathbb Z/\ell\mathbb Z)^{2n}.
		\]
		Hence the $g_r$ are pairwise distinct, and
		$\chi_{n+i}=\chi_i^{-1}$ for $1\le i\le n$.
		For each $i$, from the relation
		$L_i^{-1}E_i=q_{ii}E_iL_i^{-1}$,
		\begin{align*}
			x_ix_{n+i}-q_{ii}^{-1}x_{n+i}x_i
			&=(E_iF_i-F_iE_i)L_i^{-1}\\
			&=\frac{q_{ii}}{q_{ii}-1}
			(K_iL_i^{-1}-1)\\
			&=-\frac{q_{ii}}{q_{ii}-1}
			(1-g_ig_{n+i}).
		\end{align*}
		Hence
		\[
		\lambda_{i,n+i}=-\frac{q_{ii}}{q_{ii}-1},
		\qquad
		\lambda_{i,n+i}(1-g_ig_{n+i})\ne0.
		\]
		The linking pairs $\{i,n+i\}$ cover all $2n$ simple
		indices. The assertion follows by Theorem~\ref{thm:cartan}.
	\end{proof}

	\begin{corollary}\label{cor:cartan-parameters}
		Set
		\begin{align*}
			\mathcal R
			&=\{\beta\in\Phi_J^+,\ J\in\mathcal X:
			g_\beta^{N_J}\ne1,\
			\chi_\beta^{N_J}=\eps\},\\
			\mathcal L
			&=\{(i,j):i<j,\ i\not\sim j,\
			g_ig_j\ne1,\ \chi_i\chi_j=\eps\}.
		\end{align*}
		Let $\mathscr P$ be the affine space with coordinates
		$\mu_\beta$, $\beta\in\mathcal R$, and
		$\lambda_{ij}$, $(i,j)\in\mathcal L$.
		All other parameters are zero by the restrictions
		above.
		
		Let $\mathscr P_{\mathrm{cg}}$ be the set of parameters
		for which $A^*$ is generated by its coradical.
		Then $\mathscr P_{\mathrm{cg}}$ is nonempty if and only if
		\[
		\bigcup_{\beta\in\mathcal R}\supp\beta
		\ \cup\
		\bigcup_{(i,j)\in\mathcal L}\{i,j\}
		=
		\{1,\ldots,\theta\}.
		\]
		If $\mathscr P_{\mathrm{cg}}$ is nonempty, then
		$\mathscr P_{\mathrm{cg}}$ is dense open in $\mathscr P$.
	\end{corollary}
	
	\begin{proof}
		By Lemma~\ref{lem:cartan-root-parameter}, the union of the
		supports of roots with $u_\beta(\mu)\ne0$ equals the union
		of the supports of roots with $\mu_\beta\ne0$.
		For $(i,j)\in\mathcal L$, we have
		$1-g_ig_j\ne0$, and hence $\lambda_{ij}(1-g_ig_j)\ne0$ if and only if $\lambda_{ij}\ne0$.
		For each $i$, let $\mathscr P_i$ be the subset of
		$\mathscr P$ on which at least one of the following
		conditions holds:
		\begin{itemize}
			\item $\mu_\beta\ne0$, for some $\beta\in\mathcal R$ with $i\in\supp\beta$;
			\item $\lambda_{rs}\ne0$, for some $(r,s)\in\mathcal L$ with $i\in\{r,s\}$;
		\end{itemize}
		By Theorem~\ref{thm:cartan}, 
		\[
		\mathscr P_{\mathrm{cg}}
		=\bigcap_{i=1}^{\theta}\mathscr P_i.
		\]
		Each $\mathscr P_i$ is open in $\mathscr P$ and is
		nonempty precisely when $i$ belongs to
		\[
		\bigcup_{\beta\in\mathcal R}\supp\beta
		\ \cup\
		\bigcup_{(r,s)\in\mathcal L}\{r,s\}.
		\]
		Since $\mathscr P$ is an affine space, $\mathscr P$
		is irreducible. Therefore, if every $\mathscr P_i$
		is nonempty, the finite intersection of them is nonempty
		and dense open. If some $\mathscr P_i$ is empty,
		then $\mathscr P_{\mathrm{cg}}$ is empty.
	\end{proof}
	
	\begin{remark}
		For  the classification result in \cite[Theorem 0.1]{ASCartan}, the assumption that
		the pairs $(g_i,\chi_i)$ are distinct is satisfied.
		Indeed, suppose that
		$(g_i,\chi_i)=(g_j,\chi_j)$ for some $i\ne j$.
		Then $q_{ij}=q_{ji}=q_{ii}$, 
		and the condition $q_{ij}q_{ji}=q_{ii}^{a_{ij}}$ implies $q_{ii}^{2-a_{ij}}=1.$
		For a finite Cartan matrix,
		$a_{ij}\in\{0,-1,-2,-3\}$ when $i\ne j$.
		Thus the order of $q_{ii}$ divides one of
		$2,3,4,5$. On the other hand, the order of
		$q_{ii}=\chi_i(g_i)$ divides $|\Gamma|$.
		This contradicts the assumption that every prime
		divisor of $|\Gamma|$ is greater than $7$.
	\end{remark}

	\section{Coradicals of the dual}\label{sec:blocks}
	Let $k$ be algebraically closed of characteristic
	zero, and $A=A_D$ a lifting of a quantum linear
	space described in Subsection~\ref{sec:qls}.
	Beattie~\cite[Section 3]{Beattie} studied the coalgebra
	structure of $A^*$ under certain assumptions. We determine the coradicals of the resulting summands
	when the generators have a common character of order two,
	and examine when PBW coordinate functionals belong to
	the coradical and to $B(A^*)$.
	
	\subsection{Central idempotents and block decomposition}
	
	\begin{definition}\cite[Section 3]{Beattie}
		Let $\mathcal Z=\{x_1^{e_1}\cdots x_n^{e_n}\mid 0\le e_i<N_i\}$.
		For $\eta\in\widehat\Gamma$ and $z\in\mathcal Z$, define
		$\omega(\eta,z)\in A^*$ by
		\begin{equation}\label{eq:omega}
			\omega(\eta,z)(hz')=\eta(h)\delta_{z,z'}, \quad h\in\Gamma,\ z'\in\mathcal Z.
		\end{equation}
	\end{definition}
	
	Clearly, the set $\{\omega(\eta,z)\mid \eta\in\widehat\Gamma, z\in\mathcal Z \}$ forms a basis of $A^*$. We determine $G(A^*)$, which has been studied in \cite[Proposition 3.2]{Beattie} with certain assumptions.
	\begin{proposition}\label{prop:qls-grouplikes}
		The grouplike elements of $A^*$ are
		\[
		G(A^*)=
		\left\{
		\omega(\eta,1):
		\begin{array}{ll}
			\eta\in\widehat\Gamma,\\
			a_i(1-\eta(g_i^{N_i}))=0&\text{for all }i,\\
			c_{ij}(1-\eta(g_ig_j))=0&\text{for all }i<j.
		\end{array}
		\right\}.
		\]
	\end{proposition}
	\begin{proof}
		Recall that the elements in $G(A^*)$ correspond to the algebra characters of $A$. Let $\rho:A\to k$ be an algebra character. Since
		$q_{ii}=\chi_i(g_i)\ne1$, $\chi_i\neq\eps$ and there is $h\in\Gamma$ such that
		$\chi_i(h)\ne1$. From $\rho(hx_ih^{-1}-\chi_i(h)x_i)=0$, we have $\rho(x_i)=0$.
		Thus $\rho$ is completely determined by its restriction
		$\eta=\rho|_\Gamma\in\widehat\Gamma$.
		
		Applying $\rho$ to the power relations $x_i^{N_i}=a_i(1-g_i^{N_i})$, one obtains 
		\[
		a_i(1-\eta(g_i^{N_i}))=0, \qquad 1\le i\le n.
		\]
		Similarly, the linking relations $x_ix_j-q_{ij}x_jx_i=c_{ij}(1-g_ig_j)$ force 
		\[
		c_{ij}(1-\eta(g_ig_j))=0, \qquad 1\le i<j\le n.
		\]

		Conversely, let $\eta\in\widehat\Gamma$ satisfy the conditions listed in the proposition. Recall that 
		\[
		A_D=(T(V)\#k\Gamma)/(r-\kappa_D(r):r\in\mathcal W_{\mathrm{qls}})
		\]
		for some compatible datum $D$. Let $\varepsilon_T:T(V)\to k$ be the algebra augmentation. Define the map 
		\[
		\widetilde\rho:T(V)\# k\Gamma\longrightarrow k,
		\qquad
		\widetilde\rho(t\#h)=\varepsilon_T(t)\eta(h),
		\]
		for $t\in T(V)$ and $h\in\Gamma$.
		
		Since $\varepsilon_T(h\cdot t)=\varepsilon_T(t)$
		for all $h\in\Gamma$ and $t\in T(V)$, it follows that
		\[
		\begin{aligned}
			\widetilde\rho((t\#h)(s\#l))
			&=\varepsilon_T(t(h\cdot s))\eta(hl)\\
			&=\varepsilon_T(t)\varepsilon_T(s)\eta(h)\eta(l)\\
			&=\widetilde\rho(t\#h)\widetilde\rho(s\#l).
		\end{aligned}
		\]
		Namely, $\widetilde\rho$ is an algebra map. One can easily check that $\widetilde\rho(r-\kappa_D(r))=0$ for every $r\in\mathcal W_{\mathrm{qls}}$. Thus $\widetilde\rho$ induces an algebra map
		$\rho$ from $A$ to $k$. By construction, $\rho(h)=\eta(h)$ for $h\in\Gamma$, $\rho(x_i)=0$ and furthermore, $\rho=\omega(\eta,1)$.
	\end{proof}
	
	Put $\Omega=\{\eta\in\widehat\Gamma:\omega(\eta,1)\in G(A^*)\}$. We study the block decomposition of the coalgebra structure of $A^*$. Set
	\[
	\Theta=\langle\chi_1,\ldots,\chi_n\rangle
	\subseteq\widehat\Gamma,
	\qquad
	\Lambda=\bigcap_{i=1}^n\ker\chi_i.
	\]
	Since $\Gamma$ is finite abelian, by duality, one has
	\[
	\Theta=\{\theta\in\widehat\Gamma:
	\theta|_\Lambda=\varepsilon\}
	=\widehat{\Gamma/\Lambda},
	\]
	where $\widehat{\Gamma/\Lambda}$ is viewed as a subgroup
	of $\widehat\Gamma$. The restriction map $\widehat\Gamma\longrightarrow\widehat\Lambda$ given by $\eta\longmapsto\eta|_\Lambda$
	for $\eta\in \widehat\Gamma$ is surjective with kernel $\Theta$.
	
	For $\xi\in\widehat\Lambda$, set
	\[
	e_\xi=\frac1{|\Lambda|}
	\sum_{h\in\Lambda}\xi(h)^{-1}h,
	\qquad
	A_\xi=Ae_\xi,
	\qquad
	C_\xi=(A_\xi)^*.
	\]
	$C_\xi$ can be regarded as a subspace of $A^*$ through the following 
	\[
	(A_\xi)^*\longrightarrow A^*,
	\qquad
	f\longmapsto\bigl(a\longmapsto f(ae_\xi)\bigr).
	\]
	
	\begin{proposition}\label{prop:qls-central-blocks}
		The elements $e_\xi$ with $\xi\in\widehat\Lambda$, are pairwise orthogonal central idempotents and satisfy 
		$\sum_{\xi\in\widehat\Lambda}e_\xi=1$. Consequently,
		\begin{equation}\label{eq:blocks}
			A=\prod_{\xi\in\widehat\Lambda}A_\xi,
			\qquad
			A^*=\bigoplus_{\xi\in\widehat\Lambda}C_\xi.
		\end{equation}
		
		For $\eta\in\widehat\Gamma$ with $\eta|_\Lambda=\xi$, the set $\{\omega(\eta\theta,z)\mid\theta\in\Theta,\ z\in\mathcal Z\}$
		forms a basis of $C_\xi$. In particular,
		\[
		\dim C_\xi=|\Theta|\prod_iN_i.
		\]
		
		Moreover, $(A^*)_0=\bigoplus_{\xi\in\widehat\Lambda}\corad C_\xi$, where
		\begin{equation}\label{eq:block-coradical}
			\corad C_\xi=J(A_\xi)^\perp\cong(A_\xi/J(A_\xi))^*.
		\end{equation}
		Here $J(A_\xi)^\perp$ denotes the annihilator
		of $J(A_\xi)$ in $C_\xi$.
	\end{proposition}
	
	\begin{proof}
		Since $\Gamma$ is abelian and $\chi_i(h)=1$ for
		every $h\in\Lambda$, each element of $\Lambda$
		commutes with the generators of $A$ and hence $k\Lambda$ is central in $A$. Furthermore, $e_\xi e_{\xi'}=\delta_{\xi,\xi'}e_\xi$ and $\sum_{\xi\in\widehat\Lambda}e_\xi=1$. Consequently, the elements $e_\xi$ are pairwise orthogonal central
		idempotents, and \eqref{eq:blocks} follows.
		
		Let $\mathcal T$ be a set of representatives of
		$\Gamma/\Lambda$. For every $h\in\Gamma$, $h$ can be expressed uniquely as $h=t\ell$ with $t\in\mathcal T$ and
		$\ell\in\Lambda$. Since
		$\ell e_\xi=\xi(\ell)e_\xi$, one has $he_\xi=\xi(\ell)te_\xi$. From the basis of $A$ and the fact that
		$e_\xi\ne0$, $\{te_\xi z:t\in\mathcal T,\ z\in\mathcal Z\}$ is a basis of $A_\xi$.

		For $\eta|_\Lambda=\xi$, $\theta\in\Theta$ and
		$z,z'\in\mathcal Z$, we have
		\[
		\omega(\eta\theta,z)(te_\xi z')
		=\frac1{|\Lambda|}\sum_{h\in\Lambda}\xi(h)^{-1}(\eta\theta)(th)\delta_{z,z'}
		=\eta(t)\theta(t)\delta_{z,z'}.
		\]
		For each fixed $z\in\mathcal Z$, the matrix of the
		evaluations of $\omega(\eta\theta,z)$ on $te_\xi z$,
		with rows indexed by $t\in\mathcal T$ and columns
		by $\theta\in\Theta$, is $\bigl(\eta(t)\theta(t)\bigr)_{t,\theta}$.
		Since $\eta(t)\ne0$, this matrix is obtained from
		the character table of $\Gamma/\Lambda$ by multiplying
		each row by $\eta(t)$. Hence it is invertible.
		Consequently, the functionals $\omega(\eta\theta,z)$
		form a basis of $C_\xi$, and $\dim C_\xi=|\Theta||\mathcal Z|=|\Theta|\prod_iN_i$.
	\end{proof}
	
	\begin{remark}
		The idempotents $e_\xi$ need not be primitive
		in the center of $A$. Under the	hypothesis that $c_{ij}\ne0$ implies
		$\chi_i=\chi_j$, the decomposition
		\eqref{eq:blocks} is \cite[Proposition 3.2(3)]{Beattie}.
	\end{remark}
	
	The restriction map induces an isomorphism
	$\widehat\Gamma/\Theta\simeq\widehat\Lambda$.
	For $\eta\in\widehat\Gamma$, write
	\[
	\bar\eta=\eta\Theta,\qquad
	C(\bar\eta)=C_{\eta|_\Lambda}.
	\]

	\subsection{Clifford algebras and coradicals}
	
	Assume throughout this subsection that $n\ge1$ and
	\[
	\chi_1=\cdots=\chi_n=\chi,\qquad
	\chi^2=\eps,\qquad \chi(g_i)=-1.
	\]
	Then $N_i=2$, $\Theta=\langle\chi\rangle$ has order two,
	$\Lambda=\ker\chi$ and $\Gamma/\Lambda\cong C_2$.
	For a one-parameter family with $n=2$ and $\Gamma=C_8$,
	Beattie computed the coradicals and comatrix bases
	of $A^*$ \cite[Example 4.6]{Beattie}.
	We determine the coradicals of the blocks $C_\xi$
	for arbitrary $n$.
	
	\begin{definition}
		Let $U$ be a finite-dimensional vector space over $k$
		and let $Q$ be a symmetric bilinear form on $U$.
		The Clifford algebra of $(U,Q)$ is
		\[
		\Cl(U,Q)=T(U)/(uv+vu-2Q(u,v)1:u,v\in U).
		\]
		Let $\langle1\rangle$ denote the one-dimensional form
		whose value on a basis vector is one.
	\end{definition}
	
	For $\xi\in\widehat\Lambda$, choose an extension
	$\eta\in\widehat\Gamma$ of $\xi$. Let $Q_\xi$ be the
	symmetric bilinear form on $k^n$ whose matrix in the
	standard basis is given by
	\begin{equation}\label{eq:Clifford-Q}
		(Q_\xi)_{ii}=a_i(1-\eta(g_i^2)),\qquad
		(Q_\xi)_{ij}=(Q_\xi)_{ji}
		=\frac{c_{ij}}2(1-\eta(g_ig_j)),\quad i<j.
	\end{equation}
	Since $g_i^2,g_ig_j\in\Lambda$, the form $Q_\xi$
	is independent of the choice of $\eta$.
	
	For each $\xi\in\widehat\Lambda$, put
	\[
	\rho_\xi=\rank Q_\xi,\qquad X_i=x_ie_\xi.
	\]
	
	\begin{theorem}\label{thm:clifford}
		For each $\xi\in\widehat\Lambda$, $A_\xi\cong\Cl(k^{n+1},Q_\xi\oplus\langle1\rangle)$.
	\end{theorem}
	
	\begin{proof}
		Choose $t\in\Gamma\setminus\Lambda$ and an extension
		$\eta$ of $\xi$. Set $T=\eta(t)^{-1}te_\xi$.
		Since $t^2\in\Lambda$ and $\eta|_\Lambda=\xi$, we have
		\[
		T^2=\eta(t)^{-2}\xi(t^2)e_\xi=e_\xi.
		\]
		Since $\chi(t)=-1$, from the defining relations of $A$,
		\[
		TX_i=-X_iT,\qquad
		X_iX_j+X_jX_i=2(Q_\xi)_{ij}e_\xi.
		\]
		Thus $T,X_1,\ldots,X_n$ satisfy the defining
		relations of $\Cl(k^{n+1},Q_\xi\oplus\langle1\rangle)$.
		
		Now we show that $A_{\xi}$ is generated by $T,X_1,\ldots,X_n$.
		Observe that every $h\in\Gamma$ has either the form $\ell$ or
		$t\ell$ for some $\ell\in\Lambda$. Since
		$\ell e_\xi=\xi(\ell)e_\xi$, one has
		\[
		he_\xi=\xi(\ell)e_\xi
		\quad\text{or}\quad
		he_\xi=\xi(\ell)\eta(t)T,
		\]
		respectively. A linear basis of $A$ is given by 
		$hx_1^{e_1}\cdots x_n^{e_n}$, with
		$h\in\Gamma$ and $e_i\in\{0,1\}$. Then  every element of $A_\xi$ is a linear
		combination of products of $T,X_1,\ldots,X_n$. Consequently, there is a surjective algebra map from $\Cl(k^{n+1},Q_\xi\oplus\langle1\rangle)$ 
		to $A_\xi$, 
		where $\dim \Cl(k^{n+1},Q_\xi\oplus\langle1\rangle)=2^{n+1}$. By Proposition~\ref{prop:qls-central-blocks}, $\dim A_\xi=|\Theta|\prod_iN_i=2^{n+1}$ and so the map is an isomorphism.
	\end{proof}
	
	\begin{proposition}\label{prop:clifford-radical}
		For each $\xi\in\widehat\Lambda$, let $I_\xi$ be the
		ideal of $A_\xi$ generated by
		\[
		\left\{\sum_i v_iX_i:
		(v_1,\ldots,v_n)\in\ker Q_\xi\right\}.
		\]
		Then $I_\xi=J(A_\xi)$, and its nilpotency index is
		$n-\rho_\xi+1$, with index one when $I_\xi=0$.
		
		Moreover,
		\[
		A_\xi/J(A_\xi)\cong
		\Cl\bigl((k^n/\ker Q_\xi)\oplus k,
		\overline Q_\xi\oplus\langle1\rangle\bigr),
		\]
		where $\overline Q_\xi$ is the nondegenerate symmetric
		bilinear form induced by $Q_\xi$ on $k^n/\ker Q_\xi$.
	\end{proposition}
	\begin{proof}
		By Theorem~\ref{thm:clifford}, we identify
		$A_\xi$ with
		$\Cl(k^{n+1},Q_\xi\oplus\langle1\rangle)$.
		Let $T\in A_\xi$ denote the image of the basis
		vector of the summand $\langle1\rangle$ under this
		identification.
		
		Put $m=n-\rho_\xi$. Choose a basis
		$v_1,\ldots,v_m$ of $\ker Q_\xi$ and a basis
		$w_1,\ldots,w_{\rho_\xi}$ of a complement $W$
		on which $Q_\xi$ is nondegenerate. Set
		\[
		Z_j=\sum_i(v_j)_iX_i,\qquad
		Y_l=\sum_i(w_l)_iX_i.
		\]
		From the defining relations of $\Cl(k^{n+1},Q_\xi\oplus\langle1\rangle)$, we have
		\[
		Z_j^2=0,\qquad Z_jZ_h=-Z_hZ_j,
		\]
		\[
		Z_jY_l=-Y_lZ_j,\qquad Z_jT=-TZ_j.
		\]
		In particular, $I_\xi$ is the ideal generated by
		$Z_1,\ldots,Z_m$.
		
		By the PBW theorem for the Clifford algebra
		$A_\xi$, a linear basis of
		$A_\xi$ is given by
		\[
		\{Z_{j_1}\cdots Z_{j_q}
		Y_{l_1}\cdots Y_{l_p}T^\delta\mid q\in I_{0,m}, p\in I_{0,\rho_{\xi}}\;
		j_1<\cdots<j_q,\;
		l_1<\cdots<l_p,\;
		\delta\in\{0,1\}\}
		\]
		where $I_{k,m}=\{n\in\mathbb Z\mid k\le n\le m\}$. For every $r\ge0$,
		\[
		I_\xi^r=
		\operatorname{span}_k\left\{
		Z_{j_1}\cdots Z_{j_q}
		Y_{l_1}\cdots Y_{l_p}T^\delta\mid q\ge r
		\right\}.
		\]
		If $m>0$, then $I_\xi^{m+1}=0$ and  $I_\xi^m\ne0$.
		If $m=0$, then $I_\xi=0$. The nilpotency
		index of $I_\xi$ is $m+1=n-\rho_\xi+1$.
		
		Setting $Z_1=\cdots=Z_m=0$ implies that
		\[
		A_\xi/I_\xi\cong
		\Cl\bigl((k^n/\ker Q_\xi)\oplus k,
		\overline Q_\xi\oplus\langle1\rangle\bigr).
		\]
		The bilinear form
		$\overline Q_\xi\oplus\langle1\rangle$ is
		nondegenerate and has dimension $\rho_\xi+1$.
		Observe that the Clifford algebra of a
		nondegenerate symmetric bilinear form of dimension
		$2s$ is $\Mat_{2^s}(k)$, whereas that of dimension
		$2s+1$ is
		$\Mat_{2^s}(k)\oplus\Mat_{2^s}(k)$.
		Then $A_\xi/I_\xi$ is semisimple.
		
		Since $I_\xi$ is nilpotent,
		$I_\xi\subseteq J(A_\xi)$ and hence $J(A_\xi)/I_\xi=J(A_\xi/I_\xi)=0$. As a result, $I_\xi=J(A_\xi)$.
	\end{proof}
	Write $\Mco(d)=\Mat_d(k)^*$ for the comatrix
	coalgebra of dimension $d^2$, with coproduct as in
	\eqref{eq:comatrix}.
	
	\begin{corollary}\label{cor:clifford-coradical}
		For each $\xi\in\widehat\Lambda$,
		\begin{equation}\label{eq:clifford-corad}
			\corad C_\xi\cong
			\begin{cases}
				\Mco(2^{(\rho_\xi+1)/2}),
				&\rho_\xi\text{ odd},\\[1mm]
				\Mco(2^{\rho_\xi/2})
				\oplus\Mco(2^{\rho_\xi/2}),
				&\rho_\xi\text{ even}.
			\end{cases}
		\end{equation}
		In particular,
		\[
		\dim(A^*)_0
		=\sum_{\xi\in\widehat\Lambda}2^{\rho_\xi+1}.
		\]
	\end{corollary}
	\begin{proof}
		By Proposition~\ref{prop:clifford-radical},
		\[
		A_\xi/J(A_\xi)\cong
		\Cl\bigl((k^n/\ker Q_\xi)\oplus k,
		\overline Q_\xi\oplus\langle1\rangle\bigr).
		\]
		Since the form $\overline Q_\xi\oplus\langle1\rangle$
		is nondegenerate and has dimension $\rho_\xi+1$, 
		\[
		A_\xi/J(A_\xi)\cong
		\begin{cases}
			\Mat_{2^{(\rho_\xi+1)/2}}(k),
			&\rho_\xi\text{ odd},\\
			\Mat_{2^{\rho_\xi/2}}(k)
			\oplus\Mat_{2^{\rho_\xi/2}}(k),
			&\rho_\xi\text{ even}.
		\end{cases}
		\]
		Since $\corad C_\xi\cong(A_\xi/J(A_\xi))^*$
		and $\Mco(d)=\Mat_d(k)^*$, \eqref{eq:clifford-corad} follows by taking linear duals.
		
		In both cases, $\dim\corad C_\xi=2^{\rho_\xi+1}$.
		By \eqref{eq:blocks}, $A^*=\bigoplus_{\xi\in\widehat\Lambda}C_\xi$.
		Consequently, $(A^*)_0=\bigoplus_{\xi\in\widehat\Lambda}\corad C_\xi$.
	\end{proof}
	
	\begin{remark}
		For $\xi\in\widehat\Lambda$, the coalgebra
		$C_\xi$ is cosemisimple if and only if $\rho_\xi=n$,
		and pointed if and only if $\rho_\xi=0$.
		Moreover, $C_\xi$ is a single comatrix coalgebra
		if and only if $\rho_\xi=n$ and $n$ is odd.
	\end{remark}
	
	\begin{remark}
		Since $I_\xi=J(A_\xi)$, the coradical filtration
		of $C_\xi$ is $\{(C_\xi)_r=(I_\xi^{r+1})^\perp\mid r\ge0\}$, 
		where the annihilator is taken in $C_\xi$.
		In particular, for $\eta|_\Lambda=\xi$ and
		$z\in\mathcal Z$, $\omega(\eta,z)\in\corad C_\xi$ if and only if $\omega(\eta,z)(I_\xi)=0$.
	\end{remark}

	\subsection{Coordinate functionals and the Hopf coradical}
	
	We return to arbitrary compatible data. Let
	$R=\bigoplus_E R_E$ and $I=(R)$, as in
	Theorem~\ref{thm:qls}. For $f\in A^*$, the
	condition $f\in B(A^*)$ is equivalent to
	$f(I)=0$. If $f\in C_\xi$, then
	$f\in(A^*)_0$ is equivalent to
	$f(J(A_\xi))=0$ by \eqref{eq:block-coradical}.
	
	Consider the following condition:
	\begin{equation}\tag{PC}\label{eq:individual-coordinate}
		\text{For every }i,\text{ there exists }
		\eta_i\in\widehat\Gamma\setminus\Omega
		\text{ such that }
		\omega(\eta_i,x_i)\in C(\bar\eta_i)_0.
	\end{equation}
	For $v=\sum_jt_jx_j$, we have $\omega(\eta,x_i)(v)=t_i$. 
	Thus $\omega(\eta,x_i)$ extracts the coefficient
	of $x_i$ in $v$.
	
	By Theorem~\ref{thm:qls}, (CG) is equivalent to
	$R_E=0$ for every isotypic component $E$.
	Since $R_E\subseteq J(A)$, every element of
	$(A^*)_0$ vanishes on $R_E$, which leads to the following result.
	\begin{proposition}\label{prop:pc-cg}
		If $A$ satisfies (PC), then $A^*$ is generated by the coradical.
	\end{proposition}
	
	\begin{proof}
		Suppose that $R_E\ne0$ for some nonzero
		isotypic component $E=\operatorname{span}\{x_j:j\in I_E\}$.
		Then there is a non-zero element $v=\sum_{j\in I_E}t_jx_j\in R_E$ and an index $i\in I_E$ with $t_i\ne0$.
		By Lemma~\ref{lem:qls-radical-intersections},
		$v\in J(A)$. Hence every element of
		$(A^*)_0=J(A)^\perp$ vanishes on $v$.
		
		By (PC),
		$\omega(\eta_i,x_i)\in C(\bar\eta_i)_0
		\subseteq(A^*)_0$, while $\omega(\eta_i,x_i)(v)=t_i\ne0$, a contradiction.
		Thus $R_E=0$ for every
		nonzero isotypic component $E$. By Theorem \ref{thm:qls}, $B(A^*)=A^*$.
	\end{proof}
	
	The converse is not true in general. We shall give an example satisfying (CG) but not (PC).
	
	\begin{example}\label{ex:pc-counterexample}
		Let $\Gamma=C_6=\langle g\rangle$ and
		$\chi\in\widehat\Gamma$ with $\chi(g)=-1$.
		Take $n=3$ with 
		\[
		(g_1,\chi_1)=(g,\chi),\quad
		(g_2,\chi_2)=(g^3,\chi),\quad
		(g_3,\chi_3)=(g^5,\chi).
		\]
		Then $q_{ij}=-1$ and $N_i=2$ for all $1\le i\le 3$.
		Let $D$ be  the compatible datum given as follows:
		\[
		(a_1,a_2,a_3)=(1,0,1),\qquad
		(c_{12},c_{13},c_{23})=(2,0,2).
		\]
		Then the Hopf algebra $A=A_D$ has dimension $48$.
	\end{example}

	\begin{proposition}\label{prop:pc-not-necessary}
		Let $A$ be the Hopf algebra in
		Example~\ref{ex:pc-counterexample}. Then
		$A^*$ satisfies (CG), but $A$ does not satisfy (PC).
		Moreover,
		\begin{equation}\label{eq:counter-corad}
			(A^*)_0\cong
			k^{\oplus2}\oplus\Mco(2)^{\oplus4}.
		\end{equation}
	\end{proposition}
	
	\begin{proof}
		Write $x=x_1$, $y=x_2$ and $z=x_3$.
		Since the Yetter--Drinfeld pairs $(g_i,\chi_i)$
		are distinct, the nonzero isotypic components are
		$E_x=kx$, $E_y=ky$ and $E_z=kz$.
		Recall that $J(A)\cap k\Gamma=0$.
		From the relations $x^2=1-g^2$ and
		$z^2=1-g^4$, we have $x\notin J(A)$ and $z\notin J(A)$.
		If $y\in J(A)$, then
		\[
		xy+yx=2(1-g^4)\in J(A)\cap k\Gamma=0,
		\]
		a contradiction. Consequently,
		\[
		J(A)\cap E_x=J(A)\cap E_y=J(A)\cap E_z=0.
		\]
		By Lemma~\ref{lem:qls-radical-intersections},
		$R_E=0$ for every nonzero isotypic component $E$.
		By Theorem~\ref{thm:qls}, $B(A^*)=A^*$.
		
		We now compute $(A^*)_0$. In this case,
		$\Lambda=\langle g^2\rangle$,
		$\Theta=\langle\chi\rangle$ and
		$\Omega=\{\varepsilon,\chi\}$.
		Identify $\widehat\Lambda$ with the set of cube roots
		of unity by $\kappa=\xi(g^2)$. The corresponding
		central idempotents are
		\[
		e_\kappa=\frac13(1+\kappa^{-1}g^2+\kappa^{-2}g^4),\qquad \kappa^3=1.
		\]
		Set $A_\kappa=Ae_\kappa$ and $C_\kappa=(A_\kappa)^*$. Put $a=1-\kappa$ and $b=1-\kappa^2$.
		By \eqref{eq:Clifford-Q}, the matrix of $Q_\kappa$ is
		\[
		Q_\kappa=
		\begin{pmatrix}
			a&b&0\\
			b&0&a\\
			0&a&b
		\end{pmatrix}.
		\]
		
		If $\kappa=1$, then $Q_\kappa=0$ and $\rank Q_\kappa=0$. Suppose that $\kappa\ne1$.
		Since $\kappa^2+\kappa+1=0$, we have
		\[
		ab=a+b=3,\qquad a^3+b^3=0.
		\]
		Consequently,
		\[
		\det Q_\kappa=0,\qquad \det\begin{pmatrix}a&b\\b&0\end{pmatrix}
		=-b^2\ne0.
		\]
		Moreover,
		\[
		Q_\kappa
		\begin{pmatrix}\kappa^2\\1\\\kappa\end{pmatrix}
		=0.
		\]
		Thus
		\[
		\rank Q_\kappa=2,\qquad
		\ker Q_\kappa=k(\kappa^2,1,\kappa).
		\]
		
		By Corollary~\ref{cor:clifford-coradical},
		\[
		\corad C_1\cong k^{\oplus2},
		\qquad
		\corad C_\kappa\cong
		\Mco(2)\oplus\Mco(2)
		\quad(\kappa\ne1).
		\]
		Since $\widehat\Lambda$ has three elements, exactly
		two of which are different from $1$,
		\eqref{eq:counter-corad} follows from the decomposition
		\[
		(A^*)_0=\bigoplus_{\kappa^3=1}\corad C_\kappa.
		\]
		
		Finally, let $\eta\in\widehat\Gamma\setminus\Omega$
		and put $\kappa=\eta(g^2)$. Since
		$\Omega=\{\varepsilon,\chi\}$, we have $\kappa\ne1$.
		In particular, $C(\bar\eta)=C_\kappa$.
		By Proposition~\ref{prop:clifford-radical},
		$J(A_\kappa)$ is generated by
		\[
		u_\kappa=e_\kappa(\kappa^2x+y+\kappa z).
		\]
		
		By the definition of the coordinate functionals
		$\omega(\eta,z)$ and the formula for $e_\kappa$,
		\[
		\omega(\eta,x)(e_\kappa x)
		=\frac13\sum_{r=0}^2\kappa^{-r}\eta(g^{2r})=1.
		\]
		The other evaluations are
		\[
		\omega(\eta,y)(e_\kappa y)
		=\omega(\eta,z)(e_\kappa z)=1,
		\]
		and
		\[
		\omega(\eta,x)(e_\kappa y)
		=\omega(\eta,x)(e_\kappa z)
		=\omega(\eta,y)(e_\kappa x)
		=\omega(\eta,y)(e_\kappa z)
		=\omega(\eta,z)(e_\kappa x)
		=\omega(\eta,z)(e_\kappa y)=0.
		\]
		Therefore,
		\[
		\omega(\eta,x)(u_\kappa)=\kappa^2,\qquad
		\omega(\eta,y)(u_\kappa)=1,\qquad
		\omega(\eta,z)(u_\kappa)=\kappa.
		\]
		Since $u_\kappa\in J(A_\kappa)$, none of the three
		functionals lies in $\corad C_\kappa$.
		This holds for every
		$\eta\in\widehat\Gamma\setminus\Omega$.
		Consequently, no choice of $\eta_i$ satisfies
		condition (PC).
	\end{proof}
	
	\begin{remark}\label{rem:without-chevalley}
		Let $H=A^*$ for the lifting in
		Example~\ref{ex:pc-counterexample}. Then
		$H=B(H)$ and
		\[
		\dim H_0=18<48=\dim H.
		\]
		Hence $H_0$ is not a Hopf subalgebra of $H$.
		
		Andruskiewitsch and Cuadra proposed the classification of
		Hopf algebras over a fixed Hopf coradical~\cite{AC}.
		Garc\'ia and Jury Giraldi studied this problem for the
		$8$-dimensional Hopf algebra $\mathcal K$, assuming that
		the diagram is a Nichols algebra and the infinitesimal
		braiding is indecomposable~\cite{GG}.
		Subsequent related work includes that of Hu and
		Xiong~\cite{HX}. These works provide various examples of finite-dimensional
		Hopf algebras whose coradicals are not subalgebras,
		including families of non-basic Hopf algebras
		\cite{Xiong12,XiongSmall,Xiong24,HXRadford,BagioGarciaGiraldiMarquez}.
	\end{remark}

	\section*{acknowledgements}
	This work was partially supported by the National Natural Science Foundation of China (Grant Nos. 12401041 and 12171155), and the Science and Technology Commission of Shanghai Municipality (Grant No. 22DZ2229014).

\end{document}